%% file: MonotonicRevised_second.tex
\documentclass[12pt]{article}
\usepackage{amsmath , amsthm , commath, epsf}
\usepackage{amssymb}
\usepackage{amsfonts}
\usepackage{pause}
\usepackage{graphicx}
\usepackage{url}

\include{macros}

\begin{document}

\title{Monotonicity of the Lozi Family and the Zero Entropy Locus}
\author{
         Izzet Burak Yildiz
                }
\date{}
\maketitle

\begin{abstract}
In \cite{Ish3}, Ishii and Sands show the monotonicity of the Lozi family $\mathcal L_{a,b}$ in a $\mathcal C^1$ neighborhood of $a$-axis in the $a$-$b$ parameter space. We show the monotonicity of the entropy in the vertical direction around $a=2$ and in some other directions for $1<a\leq2$. Also we give some rigorous and numerical results for the parameters at which the Lozi family has zero entropy.

\end{abstract}

\section{Introduction}
Since its discovery in 1976, the H\'enon map \cite{Henon} has been one of the most studied examples in dynamical systems. It was introduced by M. H\'enon as a simple model exhibiting chaotic motion. On the other hand, the Lozi map \cite{Lozi} which is a piecewise affine analog of the H\'enon map has been also important since it has a simpler structure but similar chaotic behavior.

The H\'enon family is defined by:
\[ H = H_{a,b} : \left( \begin{array}{ccc} x \\ y \end{array} \right) \mapsto \left( \begin{array}{ccc} 1-ax^2+by \\ x \end{array} \right) \textit{,} \hspace{4mm} a \textit{,} b \in \mathbb{R} \textit{,} \hspace{2mm} b \neq 0,
\]
while the Lozi family is defined by:
\[ \mathcal L = \mathcal L_{a,b} : \left( \begin{array}{ccc} x \\ y \end{array} \right) \mapsto \left( \begin{array}{ccc} 1-a|x|+by \\ x \end{array} \right) \textit{,} \hspace{4mm} a \textit{,} b \in \mathbb{R} \textit{,} \hspace{2mm} b \neq 0 .
\]
Thus, the quadratic term $ax^2$ in the H\'enon family is replaced by the piecewise affine term $a|x|$. This results in a considerably simpler family of maps. For instance, in \cite{Mis3} the existence of attractors is proved for a large set of parameters, while in the H\'enon family, this is only proven for $a$ near $2$ and $b\neq0$ small (see \cite{BC}).

In this article we improve some of the entropy results obtained by Ishii and Sands in \cite{Ish3} and give some partial results about the parameters at which the topological entropy of the Lozi family is zero.

The following result about monotonicity was obtained in \cite{Ish3}:
\begin{theorem}\label{a-monotonicity} For every $a_{*} > 1$ there exists $b_{*} > 0$ such that, for any fixed $b$ with $|b| < b_{*}$, the topological entropy of $\mathcal L_{a,b}$ is a non-decreasing function of $a>a_{*}$.
\end{theorem}

Our results can be summarized in the next three theorems:

\begin{theorem} \label{b-monotonicity}For any fixed $a^{*}$ in some neighborhood of $a=2$, there exist $b_{1}^{*} > 0$ and $b_{2}^{*} < 0 $ such that the topological entropy of $\mathcal L_{a,b}$ is a non-increasing function of $b$ for $0<b<b_{1}^{*}$ and a non-decreasing function of $b$ for $b_{2}^{*}<b<0$.
\end{theorem}

Note that when $a^{*} >2$, the proof of the above theorem is trivial since those parameters stay inside the maximal entropy region where the entropy is constant and equals log$2$ (see Fig.~\ref{loziparameterspace1}). So, the non-trivial part is the one-sided neighborhood of $2$, $a^{*} \leq 2$. \\

Let us define $\mathbb{R}^2_{>1^+} = \{(a,b) \in \mathbb{R}^2 \thinspace | \thinspace a>1+|b| \}$.

\begin{Theorem}\label{angle monotonicity}
For every $1<a\leq2$ there exist $N_a^1, N_a^2\in\mathbb{R^+}$ and two lines $\gamma_{1,2}:(-\delta_{1,2},\delta_{1,2})\to\mathbb{R}_{>1^+}^2$, $\delta_{1,2}>0$, given by $\gamma_1(t)=(a+N_a^1t, -t)$ and $\gamma_2(t)=(a+N_a^2t, t)$ such that the topological entropy of $\sL_{\gamma_1(t)}$ and $\sL_{\gamma_2(t)}$ is a non-decreasing function of $t$.
\end{Theorem}

\begin{theorem} \label{zero entropy} In a small neighborhood of the parameters $a=1$ and $b=0.5$, topological entropy of $\mathcal L_{a,b}$, $h_{top}(\mathcal L_{a,b})$, is zero.
\end{theorem}

\emph{Remark:} The proof of this last result can be extended to other parameters as well. But applying the method becomes difficult especially when $b$ is close to $1$. So we give some numerical results for such parameters and obtain a picture (see Fig.~\ref{zero parameters}) for the zero entropy locus $H_{0}=\{(a,b)|\thinspace h_{top}(\mathcal L_{a,b})=0\}$ when $a>0$ and $b>0$.

\paragraph{Outline}
The remainder of this article is organized as follows.
Section~\ref{Pruning Theory} gives an introduction to the Pruning Theory and some results by Ishii and Sands that we are going to use.
Our monotonicity results are proved in Section~\ref{results1}. Then, Section~\ref{Extension} extends these results. Section~\ref{results2} describes the results about the zero entropy locus.

\section{Pruning Theory}\label{Pruning Theory}
The Pruning Theory was suggested by Cvitanovi\'{c} \cite{Cvit} as a way of obtaining symbolic dynamics for the H\'enon map. Certain conjectures were formulated which still remain unproved. Motivated by this, and following suggestions of J.~Milnor, Ishii \cite{Ish1}, \cite{Ish2} provided an analogous Pruning Theory for hyperbolic Lozi maps (i.e., those satisfying $a>1+|b|$) and proved an appropriate "Pruning Conjecture" which yielded a good symbolic description of the bounded orbits of hyperbolic Lozi maps. \\

Let us recall the basic elements of this Pruning Theory: \\

Let $\Sigma$ denote the symbol space $\{-1,+1\}^{\mathbb{Z}}$ with product topology. Define the shift map $\sigma : \Sigma \rightarrow \Sigma$ which is a continuous map given as $\sigma(\ldots \varepsilon_{-2},\varepsilon_{-1}\cdot\varepsilon_{0},\varepsilon_{1}\ldots) = (\ldots \varepsilon_{-2},\varepsilon_{-1},\varepsilon_{0}\cdot\varepsilon_{1}\ldots)$. For any $\underline{\grve} \in \Sigma$ we call $\ueps^{u} = (\ldots \varepsilon_{-2},\varepsilon_{-1})$ the tail of $\ueps$ and $\ueps^{s} = (\varepsilon_{0},\varepsilon_{1} \ldots )$ the head of $\ueps$. Let $C^u$ and $C^s$ be the set of all tails and heads, respectively. So $\Sigma$ may be identified with $C^u \times C^s$.

\medskip
Define $p(\ldots, \varepsilon_{-2},\varepsilon_{-1})(a,b)=1-bs_{-2}+b^2s_{-2}s_{-3}-b^3s_{-2}s_{-3}s_{-4}+\ldots$,

\medskip
\noindent where $s_n$ is defined as

\begin{equation} \label{s_n}
s_n \equiv \cfrac{1}{-a\varepsilon_n + \cfrac{b}{-a\varepsilon_{n-1} + \cfrac{b}{-a\varepsilon_{n-2} + \cfrac{b}{\ddots}}}}.
\end{equation}
\medskip
Similarly define $q(\varepsilon_{0},\varepsilon_{1} \ldots ) = r_0 - r_0 r_1 + r_0 r_1 r_2 - \ldots$,

\medskip
\noindent where $r_n$ is defined as
\begin{equation}\label{r_n}
r_n \equiv \cfrac{1}{a\varepsilon_n + \cfrac{b}{a\varepsilon_{n+1} + \cfrac{b}{a\varepsilon_{n+2} + \cfrac{b}{\ddots}}}}.
\end{equation}

Note that $p(\ueps^u)(a,b)$ and $q(\ueps^s)(a,b)$ are defined on $C^u \times \mathbb{R}^2_{>1^+}$ and $C^s \times \mathbb{R}^2_{>1^+}$, respectively. In the rest of the paper, we identify $p$ with $p\circ \tilde{\pi}_u$ and $q$ with $q\circ \tilde{\pi}_s$ where $\tilde{\pi}_u : \Sigma \times \mathbb{R}^2_{>1^+} \to C^u \times \mathbb{R}^2_{>1^+}$ is the map $(\ueps)(a,b) \to (\uepsu)(a,b)$ and $\tilde{\pi}_s : \Sigma \times \mathbb{R}^2_{>1^+} \to C^s \times \mathbb{R}^2_{>1^+}$ is the map $(\ueps)(a,b) \to (\uepss)(a,b)$. So, we consider $p$ and $q$ as functions $p,q : \Sigma \times \mathbb{R}^2_{>1^+} \to \mathbb{R}$.

For the proof of the next lemma, see lemma 4.3 and 6.1 in \cite{Ish1}.
\begin{Lemma} For fixed $\ueps \in \Sigma$, the functions $p(\ueps)$, $q(\ueps)$, $s_n(\ueps)$, $r_n(\ueps)$ $:\mathbb{R}^2_{>1^+} \to \mathbb{R}$ are real analytic in $(a,b)$. Moreover, $p$, $q$, $s_n$, $r_n$ and their partial derivatives with respect to $a$ and $b$ are continuous as functions $\Sigma \times \mathbb{R}^2_{>1^+} \to \mathbb{R}$.
\end{Lemma}

\begin{Definition}  \emph{We call}
\[\sP_{a,b} \equiv \{ \ueps \in \Sigma \thinspace | \thinspace (p-q)\fulleps(a,b)=0\}
\]
\emph{the} pruning front of $\sL_{a,b}$ \emph{and}
\[\sD_{a,b} \equiv \{\ueps \in \Sigma \thinspace | \thinspace (p-q)\fulleps(a,b)<0\}
\]
\emph{the} primary pruned region of $\sL_{a,b}$.
\emph{The pair} $(\sP_{a,b} , \sD_{a,b})$ \emph{is known as the} pruning pair of $\sL_{a,b}$.
\end{Definition}
We call $\sA_{a,b} \equiv \Sigma \setminus \bigcup_{n \in \mathbb{Z}} \sigma^n \sD_{a,b} = \{\ueps \in \Sigma \thinspace | \thinspace (p-q)(\sigma^n\ueps)(a,b) \geq 0 \thinspace \forall n \in \mathbb{Z} \}$ the \emph{admissible set}. \\

We explain by the end of this section that one can show the increase of the entropy by showing the decrease of the primary pruned region.

\begin{Definition}
\emph{The set} $\sPh_{a,b} \equiv \sP_{a,b} \cap \sA_{a,b}$ \emph{is the} admissible pruning front.
\end{Definition}
Let $K = K_{a,b}$ denote the set of all points whose forward and backward orbits remain bounded under $\mathcal L_{a,b}$.
For a point $X \in K$ we define $\pi(X)$ to be the set of sequences $\fulleps$ where
\[ \varepsilon_i \equiv \left\{ \begin{array}{ccc} +1 & \mbox{ if } & \sL^i(X)_x>0 \\
* & \mbox{ if } & \sL^i(X)_x=0 \\
-1 & \mbox{ if } & \sL^i(X)_x<0
\end{array} \right\} \textit{.}
\]
\\
Here $*$ can be both $+1$ and $-1$; and $Y_x$ is the $x$-component of $Y$. An element of $\pi(X)$ is called an itinerary of $X$. So a point $X$ can have more than one itinerary. \\

Now let us define the standard partial orders on $C^s \cup C^u$:
\begin{Definition}\emph{}
\begin{itemize}
\item[1.] Let $\uepss$ and $\underline{\delta}^s$ be two distinct elements in $C^s$. Then there exists the smallest number $i\geq0$ such that $\varepsilon_i\neq\delta_i$. We say $\uepss<_s\underline{\delta}^s$ if one of the following is satisfied:
\begin{itemize}

\item[(i)] The number of $+1$'s in $\cdot\varepsilon_0\ldots\varepsilon_{i-1}$ is even and $\varepsilon_i<\delta_i$,
\item[(ii)] The number of $+1$'s in $\cdot\varepsilon_0\ldots\varepsilon_{i-1}$ is odd and $\varepsilon_i>\delta_i$,\\
where order on the symbols is $-1<+1$.

\end{itemize}

\item[2.] Let $\uepsu$ and $\underline{\delta}^u$ be two distinct elements in $C^u$. Then there exists the largest number $i<0$ such that $\varepsilon_i\neq\delta_i$. When $b>0$ (resp. $b<0$), we say $\uepsu<_u\underline{\delta}^u$ if one of the following is satisfied:

\begin{itemize}
\item[(i)] The number of $-1$'s (resp. +1's) in $\varepsilon_{i-1}\ldots\varepsilon_0\cdot$ is even and $\varepsilon_i<\delta_i$,
\item[(ii)] The number of $-1$'s (resp. +1's) in $\varepsilon_{i-1}\ldots\varepsilon_0\cdot$ is odd and $\varepsilon_i>\delta_i$,
where order on the symbols is $-1<+1$.

\end{itemize}

\end{itemize}
\end{Definition}

\begin{figure}
\centerline{
\includegraphics[width=90 mm]{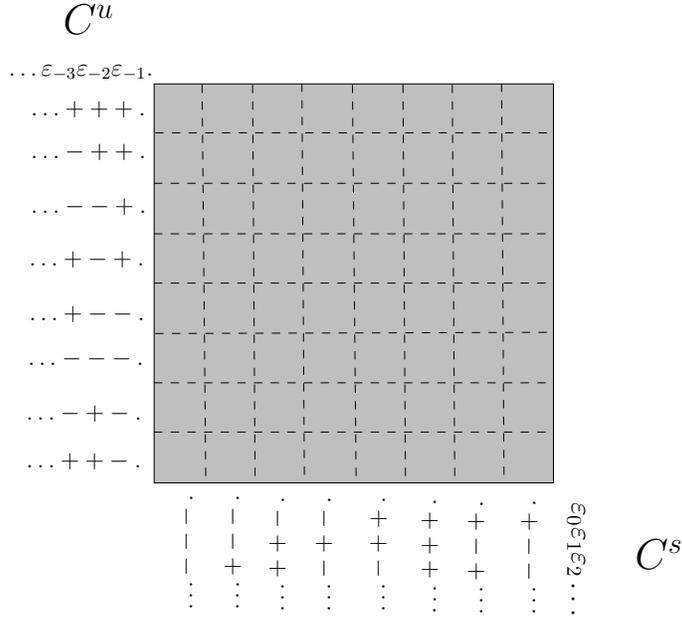}
}
\caption{Symbol space $\Sigma = \{-1,+1\}^{\mathbb{Z}}$ with the proper ordering when $b>0$. Note that $\Sigma$ can be identified with $C^u \times C^s$. To visualize $C^u \times C^s$ which is the product of two Cantor sets, we identify the endpoints of gaps (i.e., we delete the gaps) at each of the Cantor sets.}
\label{Symbol Space Intro}
\end{figure}

See Fig.\ref{Symbol Space Intro} for the case $b>0$.\\

In \cite{Ish1}, Ishii proves the following version of the Pruning Front Conjecture (PFC) which was motivated by Cvitanovi\'{c} \emph{et al} \cite{Cvit}.
\begin{Theorem}[the pruning front conjecture] Suppose that $\mathcal L_{a,b}$ satisfies $a > |b| + 1$ and
let $\ueps \in \{+1,-1\}^{\mathbb{Z}}$. Then there exists a point $X \in K_{a,b}$ such that $\ueps \in \pi(X)$ if and only if $\sigma^n \ueps$ does not lie in $\sD_{a,b}$ for all $n \in \mathbb{Z}$.
\end{Theorem}

Next, we will summarize the results of Ishii and Sands \cite{Ish3}, which prove the monotonicity of the entropy in the positive $a$-direction.\\

Recall that the tent map $T_a : \mathbb{R} \to \mathbb{R}$ is given by $T_a(x)=1-a|x|$.
\begin{Definition} \emph{An} itinerary \emph{of a point} $x \in \mathbb{R}$ \emph{under the map} $T_a$ \emph{is an element of} $i_a(x) \equiv \{\varepsilon^s \in C^s \thinspace | \thinspace \varepsilon_i T_a^i(x) \geq 0 \enspace \forall i \geq 0\}$.
\emph{We call} $\kappa(a) \equiv i_a(1)$ \emph{the} kneading invariant \emph{of} $T_a$.
\end{Definition}

\begin{Proposition}\label{projection} Suppose $1<a \leq 2$. Then $\pi_s(\sPh_{a,0}) = \kappa(a)$ where $\pi_s : \Sigma \to C^s$ is the map $\ueps \to \uepss$
\end{Proposition}

\begin{Lemma}[Stability of $\sPh$]\label{stability} Suppose $a>1+|b|$. Then for every neighborhood $U$ of $\sPh_{a,b}$ there exists a neighborhood $V$ of $(a,b)$ such that $\sPh_{\hat{a}, \hat{b}} \subset U$ for every $(\hat{a},\hat{b}) \in V$.
\end{Lemma}

\begin{Definition} We say that $(\sPh_{a,b},\sA_{a,b}) < (\sPh^{'}_{\hat{a}, \hat{b}},\sA^{'}_{\hat{a},\hat{b}})$ if $\sA_{a,b} \subset \sA^{'}_{\hat{a}, \hat{b}}$ and $\sPh^{'}_{\hat{a},\hat{b}} \cap \sA_{a,b} = \emptyset$.
\end{Definition}

The main step in the proof of the monotonicity in \cite{Ish3} is the following theorem:

\begin{Theorem}[Local Monotonicity]\label{local monotonicity} Suppose $f:(-\delta,\delta) \to \mathbb{R}_{>1+}^2$, $\delta>0$, is $C^1$ and
\[\eval[2]{\frac{d(p-q)(\ueps)f(t)}{dt}}_{t=0} > 0
\] for all $\ueps \in \sPh_{f(0)}$. Then there exists a $C^1$ neighborhood $\mathcal F$ of $f$ and a neighborhood $I$ of $0$ such that for any $C^1$ curve $g \in \mathcal F$ the map $t \in I \to (\sPh_{g(t)} , \sA_{g(t)})$ is order preserving: if $t_1, t_2 \in I$ and $t_1 < t_2$ then $(\sPh_{g(t_1)} , \sA_{g(t_1)}) < (\sPh_{g(t_2)} , \sA_{g(t_2)})$.
\end{Theorem}

It is also proven in \cite{Ish3} that if $(\sPh_{a,b} , \sA_{a,b}) < (\sPh_{\hat{a},\hat{b}} , \sA_{\hat{a},\hat{b}})$ then $h_{\textit{top}}(\mathcal L_{a,b}) \leq h_{\textit{top}}(\mathcal L_{\hat{a},\hat{b}})$. So, the above theorem can be used to show the monotonicity of the entropy. \\

In \cite{Ish3}, Ishii and Sands show that $\frac{\partial(p-q)(\ueps)(a,0)}{\partial a}>0$ for any $\ueps \in \sPh_{a,0}$. Then they use local monotonicity to prove the following:

\begin{Theorem} For every $a_*>1$ there exists $b_*>0$ such that the map $a \in (a_*,\infty) \to (\sPh_{a,b} , \sA_{a,b})$ is order preserving for all $|b|<b_*$.
\end{Theorem}
So Theorem~\ref{a-monotonicity} follows from these facts. \\

\emph{Remark:} Note the relationship between the primary pruned region, $\sD_{a,b}$, and the entropy of $\mathcal L_{a,b}$. As the primary pruned region decreases, entropy increases. This relation lies at the core of the arguments below.

\section{Results about the monotonicity of the entropy}\label{results1}
In \cite{Ish2}, Ishii mentions that although we have monotonicity in the direction given above, we do not know anything about the monotonicity in $b$ direction. We look for a solution to this question near the point $(a,b)=(2,0)$. \\

Now we want to concentrate on the point $(a,b)=(2,0)$. We will first describe the set $\sPh_{2,0}$. Using the stability of $\sPh$ this will give us some information about $\sPh_{2,b}$ for $|b|$ small. After that we will use the local monotonicity by taking $b$-derivative of $(p-q)$ to show the monotonicity in $b$-direction around $(2,0)$.

\begin{Proposition}\label{P_(2,0)}
Let $\underline{\delta}^s=(+1,-1,-1,-1\ldots)$. For $(a,b)=(2,0)$ we have $\sPh_{2,0}=\pi_s^{-1}(\underline{\delta}^s)=\{\underline{\delta}^u \cdot +1,-1,-1,-1 \ldots \thinspace | \thinspace \underline{\delta}^u \in C^u\}$ and $\sD_{2,0}=\emptyset$.
\end{Proposition}
\begin{proof}First note that by Proposition~\ref{projection}, $\pi_s(\sPh_{2,0})=\kappa(2)=(+1,-1,-1,-1\ldots)$. So, $\sPh_{2,0}\subset\pi_s^{-1}(\underline{\delta}^s)$. To prove $\pi_s^{-1}(\underline{\delta}^s)\subset\sPh_{2,0}$, we need to show that for any $\underline{\delta}^u \in C^u$ the sequence $\underline{\delta}=(\underline{\delta}^u \cdot +1,-1,-1,-1 \ldots)$ is in $\sPh_{2,0}=\sA_{2,0}\cap\sP_{2,0}$, i.e., $(p-q)(\sigma^n\underline{\delta})(2,0)\geq0$ for $n\in\mathbb{Z}$ and $(p-q)(\underline{\delta})(2,0) = 0$. Note that for an arbitrary $\ueps\in\Sigma$, $p(\uepsu)(2,0)=1$ and $r_n=\frac{1}{2\varepsilon_n}$ and $q(\uepss)(2,0)=\frac{1}{2\varepsilon_0}-\frac{1}{2^2\varepsilon_0\varepsilon_1}+\frac{1}{2^3\varepsilon_0\varepsilon_1\varepsilon_2}-\ldots +(-1)^n\frac{1}{2^{n+1}\varepsilon_0\varepsilon_1\ldots\varepsilon_n}+\ldots$. So, $q(\uepss)(2,0)$ is maximized at only $\underline{\delta}^s=(+1,-1,-1,-1\ldots)$ and its maximum value is $\Sigma_{i=1}^{i=\infty}(\frac{1}{2})^i=1$. This shows that for any $\underline{\delta}\in\pi_s^{-1}(\underline{\delta}^s)$, $(p-q)(\underline{\delta})(2,0)=0$ and $(p-q)(\sigma^n\underline{\delta})(2,0) > 0$ for $n\neq0$. This proves $\pi_s^{-1}(\underline{\delta}^s)\subset\sPh_{2,0}$ and also $\sD_{2,0}=\emptyset$.
\end{proof}

\begin{Lemma}\label{derivative positive}
\[ \eval[2]{\frac{\partial(p-q)(\ueps)(2,b)}{\partial b}}_{b=0}=\frac{1}{2\varepsilon_{-2}}
\] for $\ueps \in \sPh_{2,0}$.
\end{Lemma}

\begin{proof}
Recall that
\[p\taileps(a,b)=1-bs_{-2}+b^2s_{-2}s_{-3}-b^3s_{-2}s_{-3}s_{-4}+\ldots \textit{,}\] and
\[q\headeps(a,b) = r_0 - r_0 r_1 + r_0 r_1 r_2 - \ldots \textit{,}
\]
 where $s_n$ and $r_n$ are given by \eqref{s_n} and \eqref{r_n}.
Taking the partial derivative of $p$ with respect to $b$ we get:
\[\frac{\partial{p}}{\partial b} = -s_{-2}-bs_{-2}^{'}+2bs_{-2}s_{-3}+b^2(s_{-2}s_{-3})^{'}+ \cdots \textit{.}
\]
Since $s_n$ are analytic $\forall n\leq-2$ we obtain:
\[\eval[2]{\frac{\partial p}{\partial b}}_{b=0}=\eval[2]{-s_{-2}}_{b=0}=\frac{1}{a\varepsilon_{-2}}=\frac{1}{2\varepsilon_{-2}} \textit{.}
\]
Now for $\frac{\partial{q}}{\partial b}$; first note that for $\ueps$ such that $\uepss=(+1,-1,-1,-1\cdots)$ we have $q(\ueps)(a,b)=\frac{b}{(a+x)(b+x)}$ where $x=(a-\sqrt{a^2+4b})/2$ (see \ref{q}).\\
Since $\frac{\partial{q}}{\partial b}$ is continuous with respect to $b$; a calculation (see \ref{delq}) shows that:
\[ \eval[2]{\frac{\partial q}{\partial b}}_{b=0}=\lim_{b \to 0}{\frac{\partial q}{\partial b}}=\lim_{b \to 0}{\frac{\partial }{\partial b}}\left(\frac{b}{(a+x)(b+x)}\right)=\frac{1-\frac{2}{a}}{a(a-1)^2}.
\]
So for $a=2$ we have $\eval[2]{\frac{\partial q}{\partial b}}_{b=0}=0$.
\end{proof}
The previous lemma says that the sign of $\eval[2]{\frac{\partial(p-q)(\ueps)(2,b)}{\partial b}}_{b=0}$ depends on $\varepsilon_{-2}$.
\begin{proof}[\textbf{Proof of Theorem~\ref{b-monotonicity}}]
First let us define:
\[\mathcal{X}\equiv\{\cdots\varepsilon_{-3},+1,+1 \cdot \varepsilon_{0},\varepsilon_{1},\varepsilon_{2}\cdots\},
\]
\[\mathcal{Y}\equiv\{\cdots\varepsilon_{-3},-1,\pm1 \cdot \varepsilon_{0},\varepsilon_{1},\varepsilon_{2}\cdots\},
\]
\[\mathcal{Z}\equiv\{\cdots\varepsilon_{-3},+1,-1 \cdot \varepsilon_{0},\varepsilon_{1},\varepsilon_{2}\cdots\}.
\]

\begin{figure}
\centerline{
\includegraphics[width=90 mm]{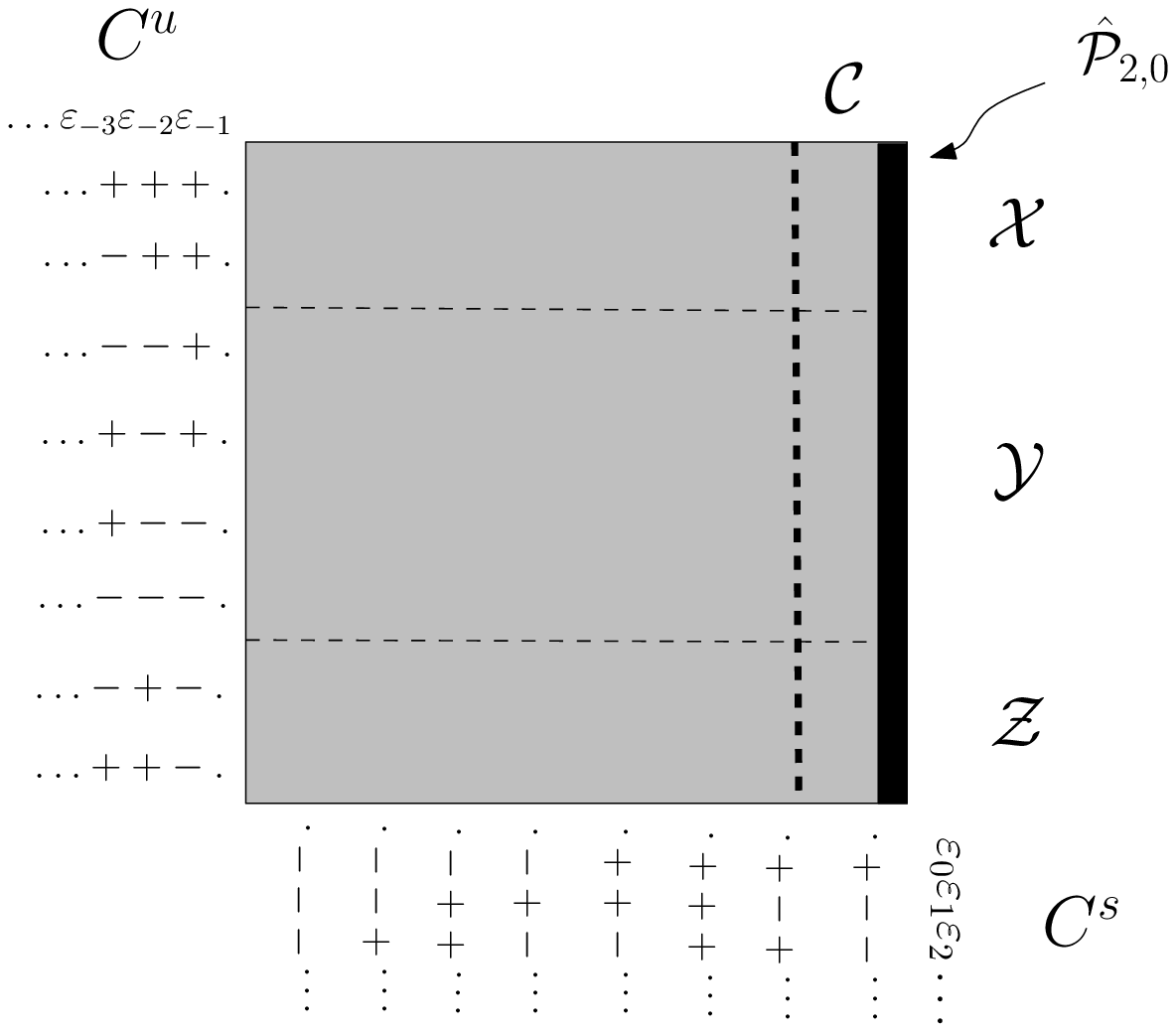}
}
\caption{Symbol space ($b>0$) and the sets $\mathcal{X}$, $\mathcal{Y}$, $\mathcal{Z}$, $\sPh_{2,0}$ and $\mathcal{C}$. Note that $\mathcal{X}$, $\mathcal{Y}$ and $\mathcal{Z}$ are the corresponding horizontal strips in the symbol space. Also note that $\sPh_{2,0}=\{\uepsu \cdot +1,-1,-1,-1 \cdots \thinspace | \thinspace \uepsu \in C^u\}$ and $\mathcal{C}$ is a neighborhood of $\sPh_{2,0}$. }
\label{Symbol Space}
\end{figure}

Also, define the curve $f(t)$ by $t \in (-\delta,+\delta) \to (2,t) \in \mathbb{R}_{>1^+}^2$ where $\delta>0$.\\
Note that we have $\sPh_{2,0}=\{\uepsu \cdot +1,-1,-1,-1 \cdots \thinspace | \thinspace \uepsu \in C^u\}$ by Proposition~\ref{P_(2,0)} and $\sD_{2,0}$ is empty (see Fig.~\ref{Symbol Space}). \\

Then by Lemma~\ref{derivative positive}, $\frac{\partial (p-q)(\ueps)(2,b)}{\partial b}$ is positive for $\ueps \in \sPh_{2,0}\cap(\mathcal{X}\cup\mathcal{Z})$ and negative for $\ueps \in \sPh_{2,0}\cap\mathcal{Y}$. \\

By continuity with respect to $\ueps$ there exists a cylinder set $\mathcal{C}$ around $\sPh_{2,0}$ such that
\[\eval[2]{\frac{\partial (p-q)(\ueps)(2,b)}{\partial b}}_{b=0}>0 \hspace{2mm} \textit{for} \hspace{2mm} \ueps \in \mathcal{C}\cap(\mathcal{X}\cup\mathcal{Z}),
 \]
 and
 \[\eval[2]{\frac{\partial (p-q)(\ueps)(2,b)}{\partial b}}_{b=0}<0 \hspace{2mm} \textit{for} \hspace{2mm} \ueps \in \mathcal{C}\cap\mathcal{Y}.
 \]
Again by continuity with respect to $b$, there exists a neighborhood $B\subset(-\delta,+\delta)$ around $0$ such that if $(2,b) \in f(B)$ we have $\frac{\partial (p-q)(\ueps)(2,b)}{\partial b}>0$ for $\ueps \in \mathcal{C}\cap(\mathcal{X}\cup\mathcal{Z})$ and $\frac{\partial (p-q)(\ueps)(2,b)}{\partial b}<0$ for $\ueps \in \mathcal{C}\cap\mathcal{Y}$.\\

Now we want to show that for $b>0$ and small, $\sPh_{2,b}\cap(\mathcal{X}\cup\mathcal{Z})$ is empty (see Fig.~\ref{pruningfronts}). \\
To do this, first observe that $\mathcal{C}$ is a neighborhood of $\sPh_{2,0}$. By stability of $\sPh$ (Lemma~\ref{stability}) there exists a neighborhood $V$ of $(2,0)$ such that $\forall (a,b) \in V \enspace \sPh_{a,b} \subset \mathcal{C}$.\\
We also know that $\frac{\partial (p-q)(\ueps)(2,b)}{\partial b}>0$ for $\ueps \in \mathcal{C}\cap(\mathcal{X}\cup\mathcal{Z})$. This means there exists a neighborhood $\tB\subset(-\delta,+\delta)$ around $0$ where $(p-q)(\ueps)$ is increasing when $b$ is increasing. This implies there exists $b_{1}^*>0$ such that for every $(2,b)$ where $0<b<b_{1}^*$ and for every $\ueps \in \mathcal{C}\cap(\mathcal{X}\cup\mathcal{Z})$ we have
\[ (p-q)(\ueps)(2,b)>(p-q)(\ueps)(2,0)\geq0.
\]
In particular, this tells us that all elements of $\sPh_{2,b}$ are in $\mathcal{C}\cap\mathcal{Y}$. But then we know that for these elements $\frac{\partial (p-q)(\ueps)(2,b)}{\partial b}<0$ and so using Theorem~\ref{local monotonicity} the entropy is non-decreasing as $b$ \emph{decreases} to $0$.\\

A similar argument applies for $b<0$ and small where it can be shown that $\sPh_{2,b}\subset \thinspace \mathcal{C}\cap(\mathcal{X}\cup\mathcal{Z})$ and that the entropy is non-decreasing as $b$ \emph{increases} to $0$.

\begin{figure}[htbp]
  \vspace{2pt}

  \centerline{\hbox{ \hspace{0.50in}
    \epsfxsize=2.5in
    \epsffile{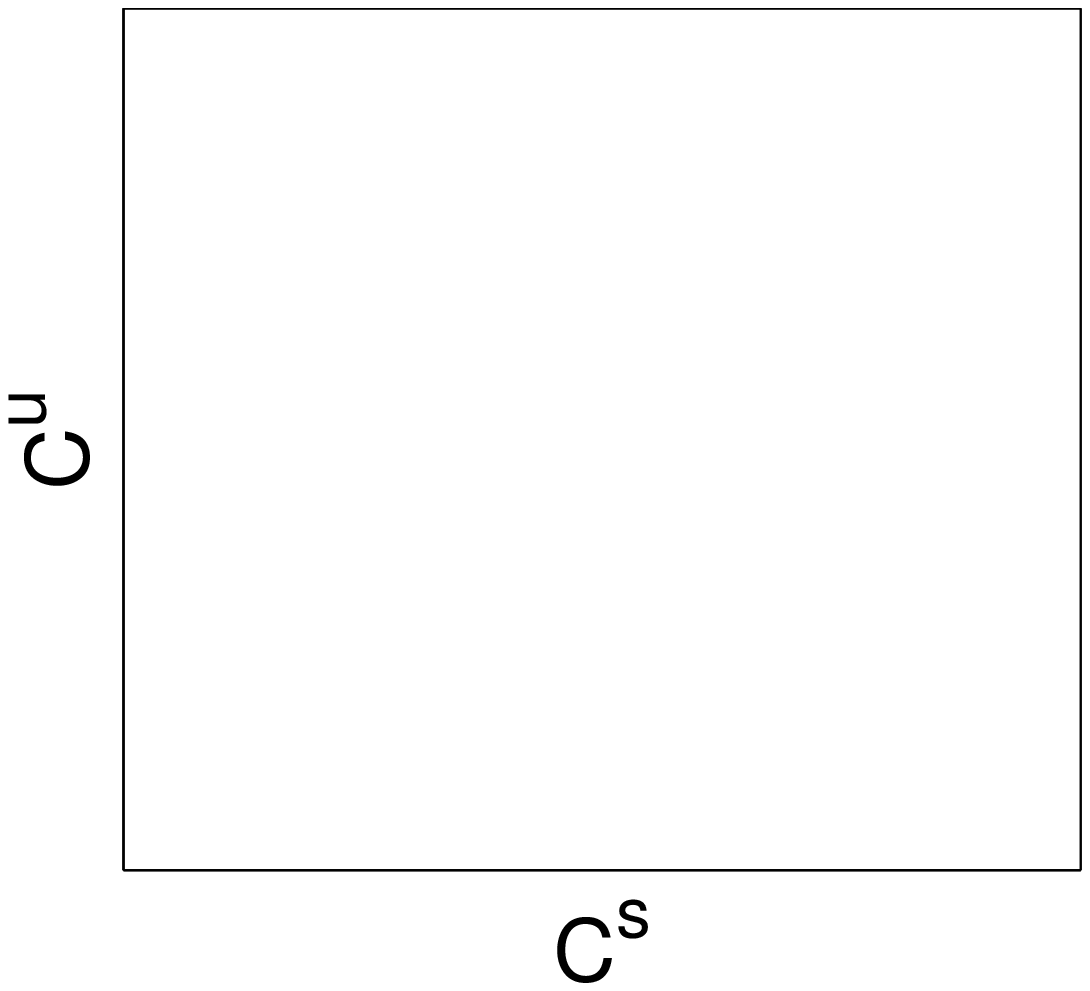}
    \hspace{0.25in}
    \epsfxsize=2.5in
    \epsffile{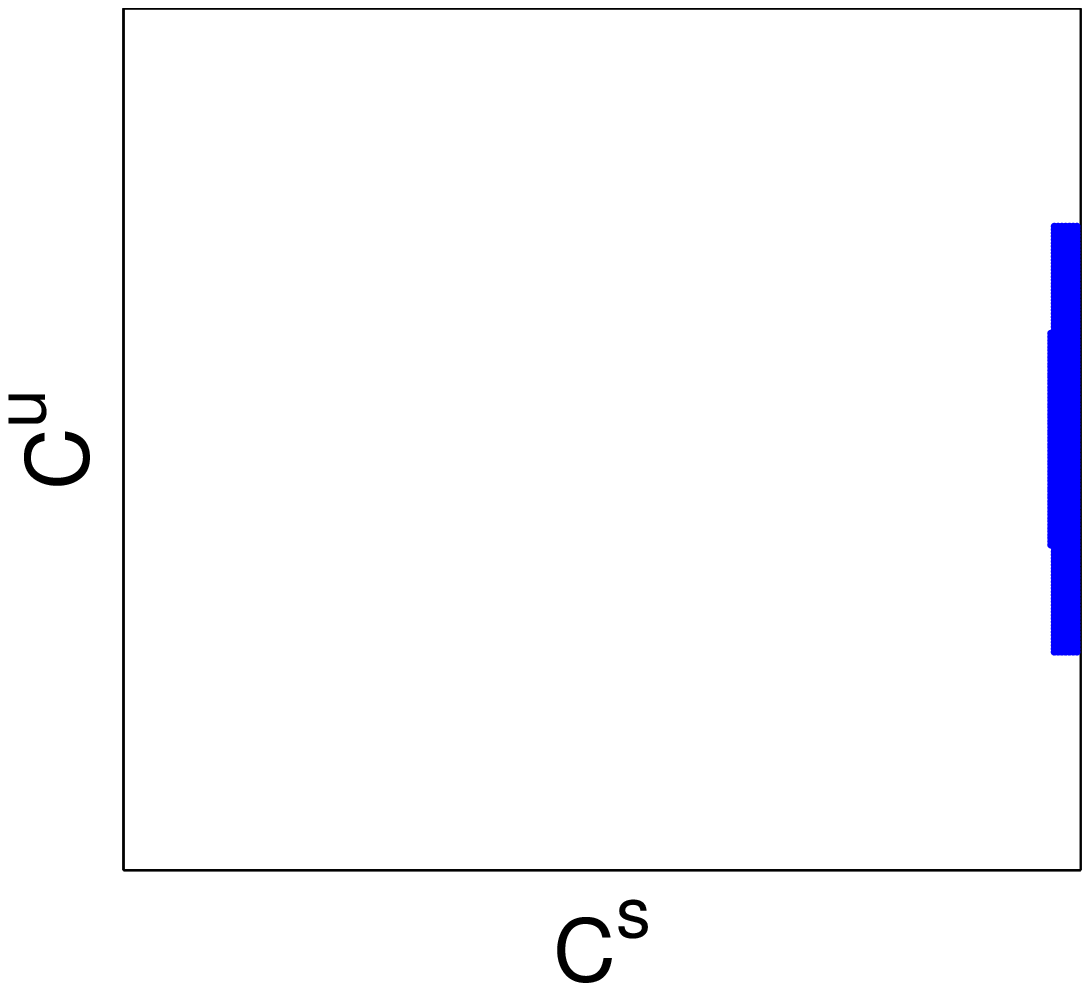}
    }
  }

  \vspace{4pt}
  \hbox{\hspace{1.3 in} (i) a=2 and b=0 \hspace{1.7in} (ii) a=2 and b=0.1}
  \vspace{10pt}

  \centerline{\hbox{ \hspace{0.50in}
    \epsfxsize=2.5in
    \epsffile{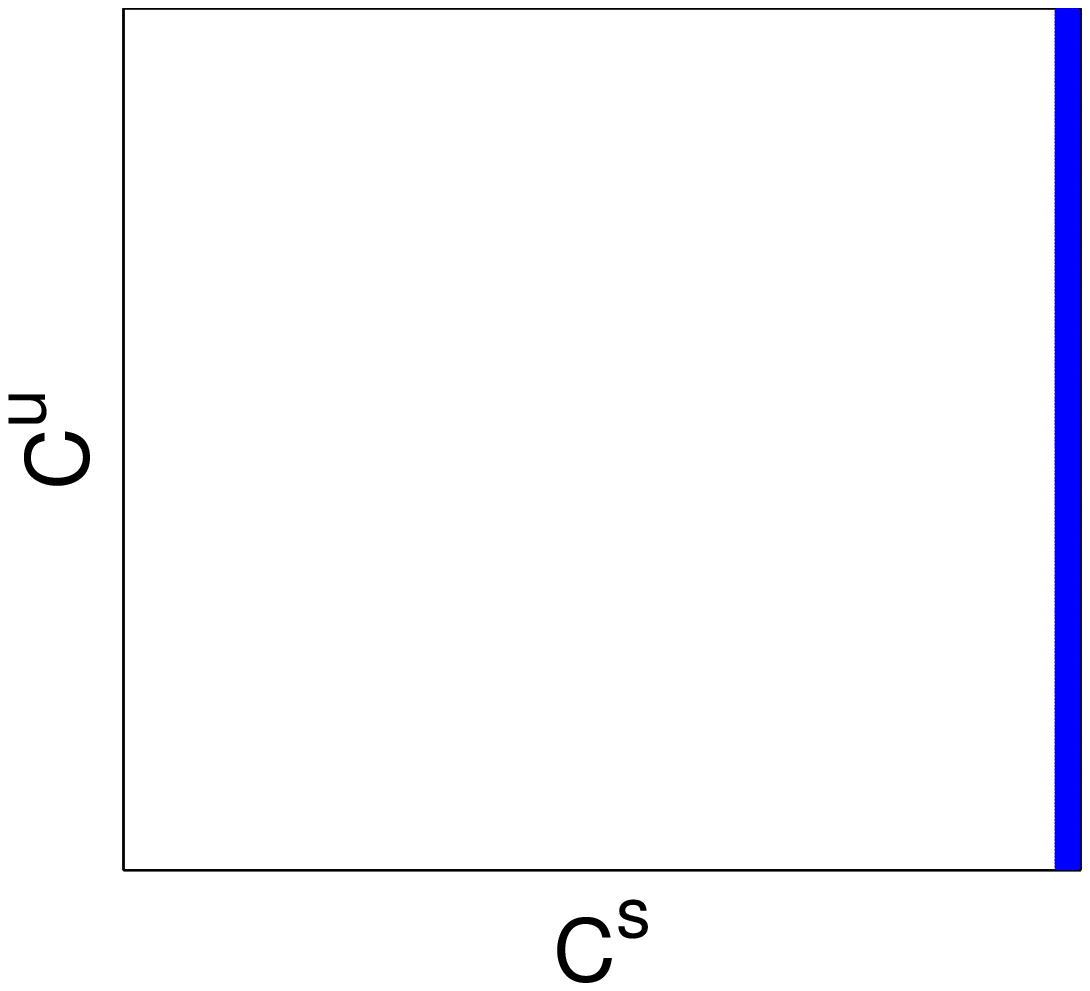}
    \hspace{0.25in}
    \epsfxsize=2.5in
    \epsffile{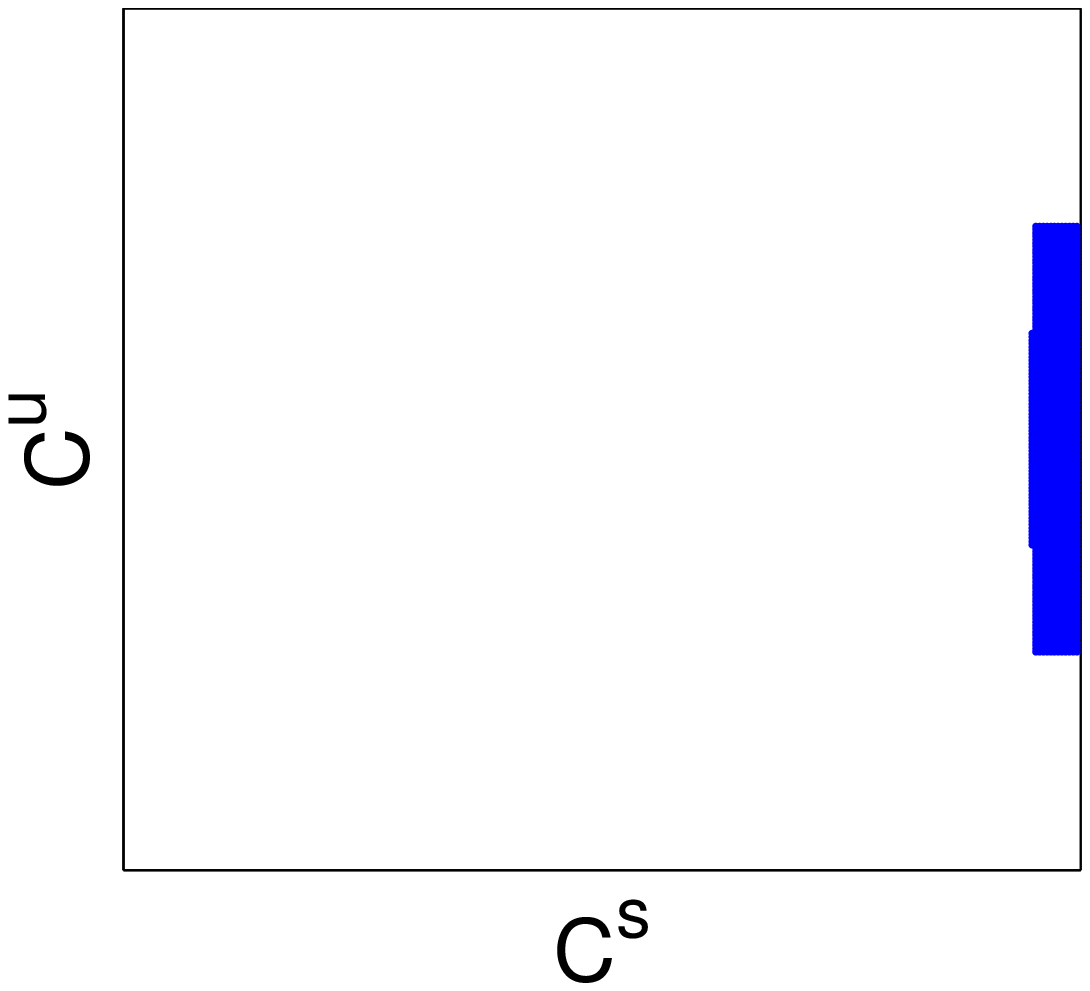}
    }
  }

  \vspace{4pt}
  \hbox{\hspace{1.3 in} (iii) a=1.95 and b=0 \hspace{1.5in} (iv) a=1.95 and b=0.1}
  \vspace{10pt}

  \caption{ This figure shows the primary pruned regions, $\sD_{a,b}$, of maps for given parameters. The $x$-axis represents $C^s$ and the $y$-axis represents $C^u$. One can expect to find some elements of $\sPh_{a,b}$ at the boundary of $\sD_{a,b}$.}

 \label{pruningfronts}

\end{figure}

\begin{figure}
\centerline{
\includegraphics[width=80mm]{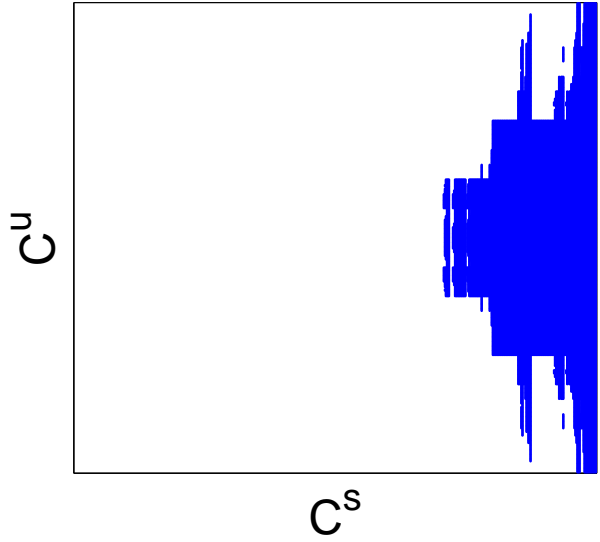}
}
\caption{Primary pruned region, $\sD_{a,b}$, for original parameters studied by Lozi: a=1.7 and b=0.5}
\label{prunedregion}
\end{figure}

\end{proof}

\section{Extension of the results to $1<a\leq2$}\label{Extension}

In this section we would like to prove some monotonicity properties for other $a$ values as well. However, we are not able to prove the monotonicity in the vertical direction because it is not possible to use local monotonicity when we move away from $a=2$. The reason behind this is the fact that for such $a$'s and small $b$, $\frac{\partial (p-q)(\ueps)(a,b)}{\partial b}$ is positive for some $\ueps\in\sPh_{a,b}$ and negative for some other $\ueps\in\sPh_{a,b}$.

\begin{figure}[hbtp]

  \centerline{\hbox{ 
    \epsfxsize=5.5in
    \epsffile{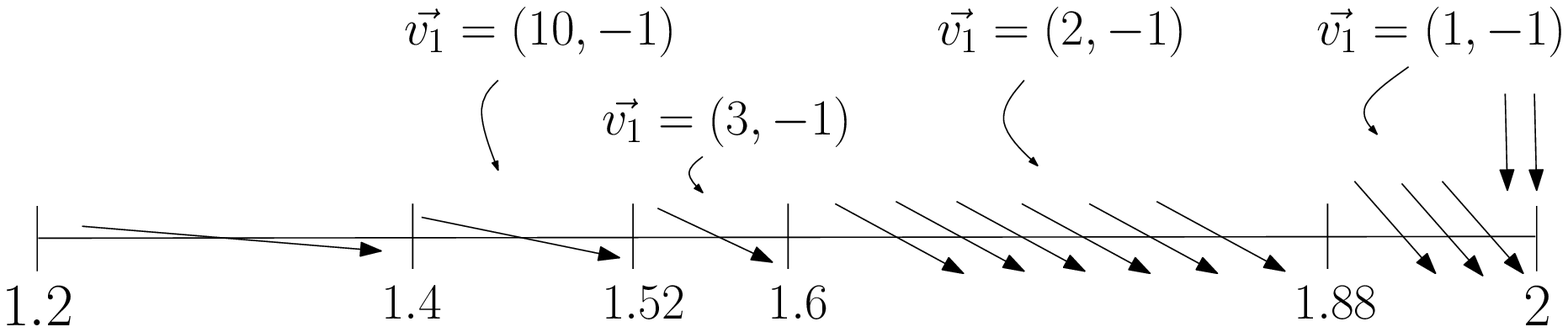}
    }
  }\vspace{0.5in}
  \centerline{\hbox{ 
    \epsfxsize=5.5in
    \epsffile{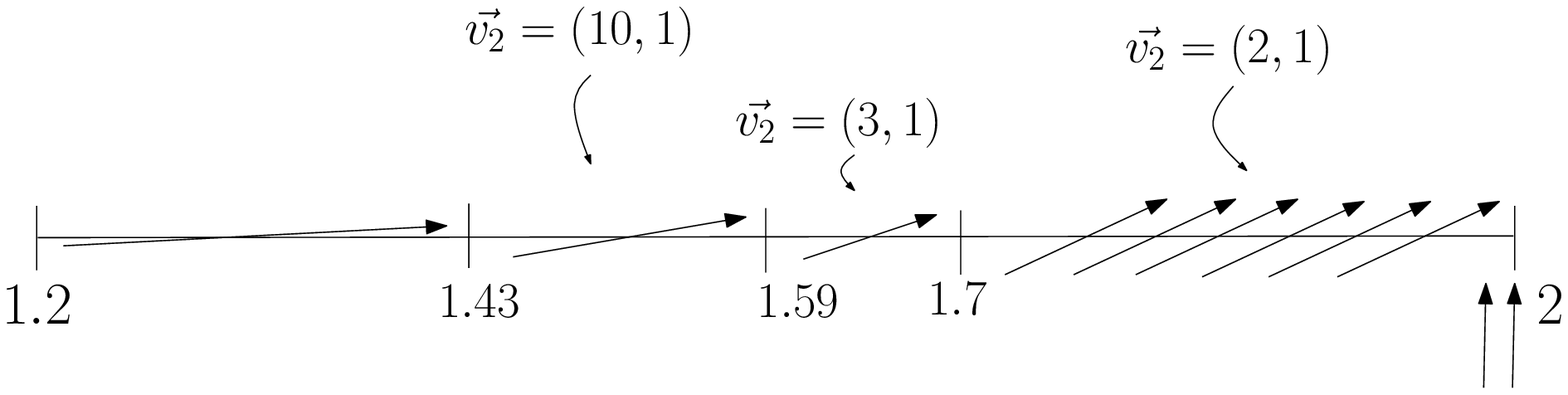}
    }
  }
  \caption{ This figure shows the approximate monotonicity results for different $N_a^1$ and $N_a^2$ values where $1.2<a\leq2$. The topological entropy is non-decreasing in the direction of arrows. }

 \label{monotonicity directions}

\end{figure}

So we prove the next best thing: Monotonicity in the direction of lines which make some angle with the $a$-axis (see Fig.~\ref{monotonicity directions}).
To prove this result we modify and use some of the computations done in \cite{Ish3}.

\begin{Lemma}(Lemma 11 in \cite{Ish3})
Suppose $1<a\leq2$ and $\uepss\in\kappa(a)$. Then
\[ \frac{a^3+2a^2-6a+2}{2a^2(a-1)} \leq \eval[2]{\frac{\partial (p-q)(\ueps)(a,b)}{\partial a}}_{(a,0)} \leq \frac{a^3+2a^2-6a+4}{2a^2(a-1)}.
\]
In particular; $\eval[1]{\frac{\partial (p-q)(\ueps)(a,b)}{\partial a}}_{(a,0)} \geq (\sqrt{2}-1)/2 >0$ if $a \geq \sqrt{2}$.
\end{Lemma}

\begin{Lemma}(Corollary 13 in \cite{Ish3})\label{a derivative positive}
Suppose $1<a\leq2$ and $\uepss\in\kappa(a)$. Then $\eval[1]{\frac{\partial (p-q)(\ueps)(a,b)}{\partial a}}_{(a,0)} > 0$.
\end{Lemma}

\begin{Lemma}(Corollary 7 and Eqn. 3.11 in \cite{Ish3})\label{equations}

Suppose $1<a\leq2$ and $\uepss\in\kappa(a)$. Then
\begin{equation}
\label{equation1}
\sum_{i=0}^{\infty}(-1)^i \frac{\varepsilon_{0}\dots\varepsilon_{i-1}}{a^i} =0,
\end{equation}
and
\begin{equation}
\label{equation2}
\sum_{j=0}^{\infty}(-1)^{i+j} \frac{\varepsilon_{0}\dots\varepsilon_{i+j}}{a^{i+j+1}} = (-1)^i\frac{\varepsilon_{0}\dots\varepsilon_{i-1}}{a^i}T_a^i(1),
\end{equation}
where we define the empty product $\varepsilon_{0} \dots \varepsilon_{-1}$ to equal 1.
\end{Lemma}

Now, we use these results and similar techniques to prove the following:
\begin{Lemma}
Suppose $1<a\leq2$ and $\uepss\in\kappa(a)$. Then
\[ \frac{1}{a\varepsilon_{-2}}- \frac{-2a^2+7a-2}{2a^3(a-1)} \leq \eval[2]{\frac{\partial (p-q)(\ueps)(a,b)}{\partial b}}_{b=0} \leq \frac{1}{a\varepsilon_{-2}}- \frac{-2a^2+7a-8}{2a^3(a-1)}.
\]
\begin{proof}
From the proof of Lemma \ref{derivative positive} we know that $\eval[1]{\frac{\partial{p}(\ueps)(a,b)}{\partial b}}_{b=0}=\frac{1}{a\varepsilon_{-2}}$. So, we need to find some upper and lower bound for $\eval[1]{\frac{\partial{q}(\ueps)(a,b)}{\partial b}}_{b=0}$.

Remember that $q(\varepsilon_{0},\varepsilon_{1} \ldots ) = r_0 - r_0 r_1 + r_0 r_1 r_2 - \ldots$. \\
Let us write $q(\varepsilon_{0},\varepsilon_{1} \ldots ) = T_{0} - T_{1} + T_{2} - T_{3} + \ldots = \sum\limits_{n=0}^{\infty} (-1)^n T_n$ where $T_n = r_0 r_1 \ldots r_n$.

Now we have the following:
\[\eval[1]{r_n}_{b=0} = \eval[1]{\frac{1}{a\varepsilon_{n}+b r_{n+1}}}_{b=0}=\frac{\varepsilon_n}{a},
\]
and
\[\eval[1]{r_n'}_{b=0} = \eval[2]{\frac{\partial (r_n)(\ueps)(a,b)}{\partial b}}_{b=0}=-\eval[1]{r_n^2 (r_{n+1}+b r_{n+1}')}_{b=0} = -\frac{\varepsilon_{n+1}}{a^3}.
\]
Taking term by term derivative of $q$, we get the following. Note that $\varepsilon_0=+1$ and $\varepsilon_1=-1$:
\[ T_0' = r_0'=-\frac{\varepsilon_1}{a^3}=\frac{1}{a^3},
\]

\[-T_1' = -(r_0'r_1+r_0r_1')=\frac{\varepsilon_1}{a^3}\frac{\varepsilon_1}{a}+\frac{\varepsilon_0}{a}\frac{\varepsilon_2}{a^3}=
\frac{1}{a^4}+\frac{\varepsilon_0\varepsilon_2}{a^4},
\]

\[ T_2' =(r_0 r_1)'r_2 + (r_0 r_1)r_2'=-\frac{\varepsilon_2}{a^5}-\frac{\varepsilon_0}{a^5}-\frac{\varepsilon_0 \varepsilon_1\varepsilon_3}{a^5},
\]

\[-T_3' = -(r_0 r_1 r_2)'r_3+(r_0 r_1 r_2)r_3'=\frac{\varepsilon_2\varepsilon_3}{a^6}
+\frac{\varepsilon_0\varepsilon_3}{a^6}+\frac{\varepsilon_0\varepsilon_1}{a^6}+\frac{\varepsilon_0\varepsilon_1\varepsilon_2\varepsilon_4}{a^6},
\]

\[ T_4' = (r_0 r_1 r_2 r_3)' r_4+(r_0 r_1 r_2 r_3)r_4'=-\frac{\varepsilon_2\varepsilon_3\varepsilon_4}{a^7}
-\frac{\varepsilon_0\varepsilon_3\varepsilon_4}{a^7}-\frac{\varepsilon_0\varepsilon_1\varepsilon_4}{a^7}
-\frac{\varepsilon_0\varepsilon_1\varepsilon_2}{a^7}-\frac{\varepsilon_0\varepsilon_1\varepsilon_2\varepsilon_3\varepsilon_5}{a^7},
\]

\[\ddots=\ddots .
\]

Note that $\eval[1]{\frac{\partial{q}(\ueps)(a,b)}{\partial b}}_{b=0} = T_0' - T_1' + T_2' - \cdots + (-1)^nT_n' + \cdots$.

{\bf Claim:} $T_n' = (r_0 r_1 \ldots r_n)' = -\frac{1}{a^{n+3}}(\sum\limits_{i=0}^{n-1}\varepsilon_0\ldots\overline{\varepsilon_i\varepsilon_{i+1}}\ldots\varepsilon_n + \varepsilon_0\varepsilon_1\ldots\overline{\varepsilon_n}\varepsilon_{n+1})$, \hspace{5mm} $n\geq2$ where $\varepsilon_0\ldots\varepsilon_{i-1}\overline{\varepsilon_i}\varepsilon_{i+1}\ldots\varepsilon_{n}$ means $\varepsilon_i$ is missing in the term. \\

\emph{Proof of the Claim:} Note that $T_n'=(T_{n-1}r_n)' = T_{n-1}'\frac{\varepsilon_n}{a} + (r_0 r_1 \ldots r_{n-1})r_n' = T_{n-1}'\frac{\varepsilon_n}{a} - \frac{\varepsilon_0\varepsilon_1\ldots\overline{\varepsilon_n}\varepsilon_{n+1}}{a^{n+3}}$. So the claim follows by induction.\\

Let us organize the terms in the following matrix form. The terms in the first column ($(-1)^nT_n'$) are equal to the sum of the terms in the corresponding row. Also note that the last row denotes the sum of the terms in the corresponding column:

\[ \left[\begin{array}{c|cccccc}
-T_1' & & & & & \frac{1}{a^4} & \frac{\varepsilon_0\varepsilon_2}{a^4} \\
T_2' & -\frac{\varepsilon_2}{a^5} & & & & -\frac{\varepsilon_0}{a^5} & -\frac{\varepsilon_0\varepsilon_1\varepsilon_3}{a^5} \\
-T_3' & \frac{\varepsilon_2\varepsilon_3}{a^6} & \frac{\varepsilon_0\varepsilon_3}{a^6} & & & \frac{\varepsilon_0\varepsilon_1}{a^6} & \frac{\varepsilon_0\varepsilon_1\varepsilon_2\varepsilon_4}{a^6} \\
T_4' & -\frac{\varepsilon_2\varepsilon_3\varepsilon_4}{a^7} & -\frac{\varepsilon_0\varepsilon_3\varepsilon_4}{a^7} & -\frac{\varepsilon_0\varepsilon_1\varepsilon_4}{a^7} & & -\frac{\varepsilon_0\varepsilon_1\varepsilon_2}{a^7} & -\frac{\varepsilon_0\varepsilon_1\varepsilon_2\varepsilon_3\varepsilon_5}{a^7} \\
\vdots & \vdots & \vdots & \vdots & \ddots & \vdots & \vdots \\
(-1)^nT_n' & -\frac{(-1)^{n}\overline{\varepsilon_0\varepsilon_1}\varepsilon_2\ldots\varepsilon_n}{a^{n+3}} &
-\frac{(-1)^{n}\varepsilon_0\overline{\varepsilon_1\varepsilon_2}\ldots\varepsilon_n}{a^{n+3}} &
\cdots &
\cdots &
-\frac{(-1)^{n}\varepsilon_0\ldots\overline{\varepsilon_{n-1}\varepsilon_{n}}}{a^{n+3}} &
-\frac{(-1)^{n}\varepsilon_0\ldots\overline{\varepsilon_{n}}\varepsilon_{n+1}}{a^{n+3}} \\
\vdots & \vdots & \vdots & \vdots & \vdots & \vdots & \vdots \\ \hline
\sum\limits_{n=1}^{\infty} (-1)^n T_n' & \star_1 & \star_2 &  \star_3 \cdots &  \star_n \cdots & S & \frac{\varepsilon_0\varepsilon_2}{a^4} + R \\

\end{array} \right]
\]

Since the series which gives the derivative of $q$ is absolutely convergent, regrouping the suitable terms together, we can write:
\[ \eval[1]{\frac{\partial{q}(\ueps)(a,b)}{\partial b}}_{b=0}=\frac{1}{a^3}+\frac{\varepsilon_0\varepsilon_2}{a^4}+\sum_{n=1}^{\infty}\star_n + S + R,
\]
where
\[ \star_1=-\frac{\overline{\varepsilon_0}\overline{\varepsilon_1}\varepsilon_2}{a^5}
+\frac{\overline{\varepsilon_0}\overline{\varepsilon_1}\varepsilon_2\varepsilon_3}{a^6}
-\frac{\overline{\varepsilon_0}\overline{\varepsilon_1}\varepsilon_2\varepsilon_3\varepsilon_4}{a^7}
+\dots,
\]
\[\star_2=\frac{\varepsilon_0\overline{\varepsilon_1}\overline{\varepsilon_2}\varepsilon_3}{a^6}
-\frac{\varepsilon_0\overline{\varepsilon_1}\overline{\varepsilon_2}\varepsilon_3\varepsilon_4}{a^7}
+\frac{\varepsilon_0\overline{\varepsilon_1}\overline{\varepsilon_2}\varepsilon_3\varepsilon_4\varepsilon_5}{a^8}
- \dots,
\]

\[\star_n=(-1)^n\frac{\varepsilon_0\ldots\overline{\varepsilon_{n-1}}\overline{\varepsilon_n}\varepsilon_{n+1}}{a^{n+4}}
+(-1)^{n+1}\frac{\varepsilon_0\ldots\overline{\varepsilon_{n-1}}\overline{\varepsilon_n}\varepsilon_{n+1}\varepsilon_{n+2}}{a^{n+5}}
+\ldots,
\]
and
\[ S=\frac{1}{a^4}-\frac{\varepsilon_0}{a^5}+\frac{\varepsilon_0\varepsilon_1}{a^6}-\ldots+(-1)^{i+1}\frac{\varepsilon_0\ldots\varepsilon_i}{a^{i+5}}+\ldots \hspace{5mm} i\geq-1,
\]
and
\[R = -\frac{\varepsilon_0\varepsilon_1\overline{\varepsilon_2}\varepsilon_3}{a^5}
+\frac{\varepsilon_0\varepsilon_1\varepsilon_2\overline{\varepsilon_3}\varepsilon_4}{a^6}
-\ldots+(-1)^{j+1}\frac{\varepsilon_0\ldots\varepsilon_{j-1}\overline{\varepsilon_j}\varepsilon_{j+1}}{a^{j+3}} \hspace{5mm} j\geq2.
\]

First let us start with observing that by the equation (\ref{equation1}):
\begin{equation}\label{S}
S=\frac{1}{a^4}(1-\frac{\varepsilon_0}{a}+\frac{\varepsilon_0\varepsilon_1}{a^2}-\frac{\varepsilon_0\varepsilon_1\varepsilon_2}{a^3}+\dots) = 0.
\end{equation}
Secondly,
\begin{equation}\label{R}
|R| \leq \frac{1}{a^5}(1+\frac{1}{a}+\frac{1}{a^2}+\dots)=\frac{1}{a^4(a-1)}.
\end{equation}
For $\sum_{n=1}^{\infty}\star_n$, by the equation (\ref{equation2}) we have:
\[ \star_1(-a^2\varepsilon_0\varepsilon_1)=(-1)^2\frac{\varepsilon_0\varepsilon_1}{a^2}T_a^2(1),
\]
\[ \star_2(-a^2\varepsilon_1\varepsilon_2)=(-1)^3\frac{\varepsilon_0\varepsilon_1\varepsilon_2}{a^3}T_a^3(1),
\]
and
\[ \star_n(-a^2\varepsilon_{n-1}\varepsilon_n)=(-1)^{n+1}\frac{\varepsilon_0\varepsilon_1\ldots\varepsilon_n}{a^{n+1}}T_a^{n+1}(1).
\]
So,
\[\sum_{n=1}^{\infty}\star_n=-\frac{1}{a^4}\sum_{i=0}^{\infty}(-1)^i\frac{\varepsilon_0\varepsilon_1\ldots\varepsilon_{i-1}}{a^i}T_a^{i+2}(1).
\]
Again by (\ref{equation1}) we also have
\[\sum_{n=1}^{\infty}\star_n=-\frac{1}{a^4}\sum_{i=0}^{\infty}(-1)^i\frac{\varepsilon_0\varepsilon_1\ldots\varepsilon_{i-1}}{a^i}(T_a^{i+2}(1)-\alpha)
\]
for any $\alpha\in\mathbb{R}$. Let $\alpha=\frac{1+T_a(1)}{2}$ and $\delta=\frac{1-T_a(1)}{2}=\frac{a}{2}$. Since $T_a^i(1)\in[T_a(1),1]$ for every $i\geq0$ we have $-\delta\leq T_a^i(1)-\alpha\leq\delta$ for every $i\geq0$. Note that by direct calculation $T_a^2(1)-\alpha=\delta(2a-3)$. This gives us:
\begin{equation}
\sum_{n=1}^{\infty}\star_n=-\frac{1}{a^4}(T_a^2(1)-\alpha -
\frac{\varepsilon_0}{a}(T_a^3(1)-\alpha)+\frac{\varepsilon_0\varepsilon_1}{a^2}(T_a^4(1)-\alpha)-\ldots) \nonumber
\end{equation}
\begin{equation}\label{sumstar1}
\leq -\frac{1}{a^4}(\delta(2a-3)-\delta(\frac{1}{a}+\frac{1}{a^2}+\frac{1}{a^3}+\ldots))=-\frac{2a^2-5a+2}{2a^3(a-1)}.
\end{equation}
Similar calculations show that
\begin{equation}\label{sumstar2}
\sum_{n=1}^{\infty}\star_n \geq -\frac{2a^2-5a+4}{2a^3(a-1)}.
\end{equation}
Now, combining (\ref{S}), (\ref{R}), (\ref{sumstar1}) and (\ref{sumstar2}) we get the desired result.
\end{proof}
\end{Lemma}

\begin{proof}[Proof of the Theorem \ref{angle monotonicity}]
By Lemma \ref{a derivative positive} we know that for any $1<a\leq2$ and $\ueps^{s}\in \kappa(a)$, $\eval[1]{\frac{\partial (p-q)(\ueps)(a,b)}{\partial a}}_{(a,0)} > 0$. Also by the previous lemma for any such $a$, $\eval[1]{\frac{\partial (p-q)(\ueps)(a,b)}{\partial b}}_{(a,0)}$ has an upper and lower bound. So, there exist $N_a^1\in\mathbb{R^+}$ such that
\[N_a^1\eval[2]{\frac{\partial (p-q)(\ueps)(a,b)}{\partial a}}_{(a,0)}-\eval[2]{\frac{\partial (p-q)(\ueps)(a,b)}{\partial b}}_{(a,0)}>0,
\]
and $N_a^2\in\mathbb{R^+}$ such that
\[N_a^2\eval[2]{\frac{\partial (p-q)(\ueps)(a,b)}{\partial a}}_{(a,0)}+\eval[2]{\frac{\partial (p-q)(\ueps)(a,b)}{\partial b}}_{(a,0)}>0.
\]
This means that the directional derivatives of $(p-q)(\ueps)(a,b)$ in the direction $\vec{v_1}=(N_a^1,-1)$ and $\vec{v_2}=(N_a^2,1)$ are both positive. So by local monotonicity theorem, result follows.
\end{proof}

Since we have explicit upper and lower bounds for both $\frac{\partial (p-q)(\ueps)(a,b)}{\partial a}$ and $\frac{\partial(p-q)(\ueps)(a,b)}{\partial b}$ we can compute the directions in which the entropy is non-decreasing (see Fig.~\ref{monotonicity directions}).

\section{Results about the zero entropy locus}\label{results2}

In this section we turn our attention to the parameters for which $h_{top}(\mathcal L_{a,b})=0$. Note that it is enough to consider the maps with $|b|\leq1$ since the maps with $|b|>1$ are, up to affine conjugacy, inverses of the maps with $|b|<1$. \\

\begin{figure}[htbp]
\centerline{
\includegraphics[width=150mm]{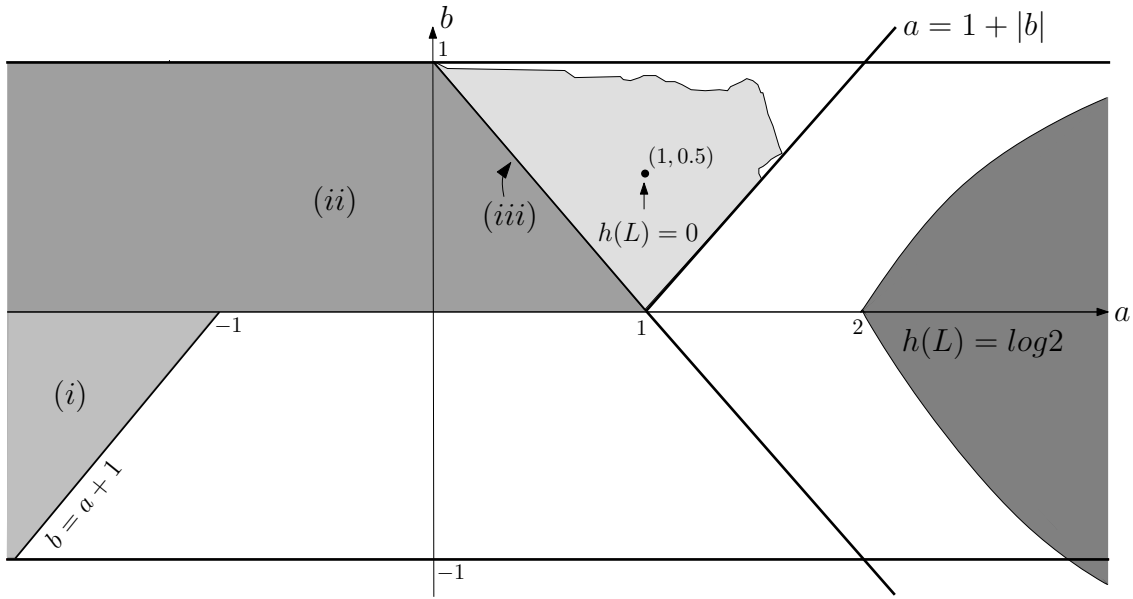}
}
\caption{The sets $(i)$, $(ii)$ and $(iii)$ where the entropy is zero (Thm.~\ref{zero entropy Ishii}) are shown above. It was also numerically observed that the light gray region with complicated boundary has zero entropy. We rigorously prove this at $(1,0.5)$. The darker gray region on the very right is the maximal entropy region where the entropy is log$2$ (see \cite{Ish1}).}
\label{loziparameterspace1}
\end{figure}

Let us first review the following theorem:
\begin{Theorem}[\cite{Ish5}]\label{zero entropy Ishii} If the Lozi map $\mathcal L_{a,b}$ satisfies either (i) $-1\leq b<0$ and $a\leq b-1$, (ii) $0<b\leq1$ and $a< 1-b$, (iii) $0<b\leq1$ and $a=1-b$, then $h_{top}(\mathcal L_{a,b})=0$ (see Fig.~\ref{loziparameterspace1} for related regions).
\end{Theorem}

\begin{proof}
 If $a\leq b-1 \leq -a$ then $\mathcal L_{a,b}$ has no fixed points. When $b<0$, $\mathcal L_{a,b}$ is orientation preserving, so by Brouwer's translation theorem \cite{Brouwer} it has an empty non-wandering set and therefore zero entropy, proving (i). Part (ii) can be investigated in two parts. When $0<b\leq1$ and $a \leq b-1$, $\mathcal L_{a,b}$ has no fixed points and no period-two points. So, one can apply Brouwer's translation theorem to $\mathcal L_{a,b}^2$ (which is orientation preserving when $b>0$) to deduce that $h_{top}(\mathcal L_{a,b})=h_{top}(\mathcal L_{a,b}^2)/2=0$. When $0<b\leq1$ and $b-1< a < 1-b$, there exists a unique saddle fixed point $p=(1/(1+a-b),1/(1+a-b))$ in the first quadrant. Also note that there is no period-two point. Now $v^s=(\lambda,1)$ where $\lambda=(-a+\sqrt{a^2+4b})/2$ is a stable direction at $p$ and $W^{s}_{+}(p)=\{p+v^{s}t \in \mathbb{R}^2 | \thinspace t>0\}$ is invariant under $\mathcal L_{a,b}$. Also $\mathbb{R}^2 \setminus (W^{s}_{+}(p) \cup \{p\} )$ is homeomorphic to $\mathbb{R}^2$ and $\mathcal L_{a,b}^2$ has no fixed points there. Since $\mathcal L_{a,b}^2$ is orientation preserving when $b>0$, $h_{top}(\mathcal L_{a,b})=h_{top}(\mathcal L_{a,b}^2)/2=0$, proving (ii). When $a=1-b$, note that there is a unique fixed point $p=(1/2a,1/2a)$ and a line segment $I$ of period-two points where $I = \{(x,y)\in \mathbb{R}^2 \thinspace | \thinspace a(x+y) =1, \thinspace 0\leq x\leq1/a \}$. Note that $p$ is the midpoint of $I$. So, one can similarly apply Brouwer's translation theorem to $\mathcal L_{a,b}^2$ on $\mathbb{R}^2 \setminus (W^{s}_{+}(p) \cup I)$ to conclude that $\mathcal L_{a,b}$ has zero entropy, proving (iii).
\end{proof}

Note that nothing much is known when $-1\leq b<0$ except the case $a\leq b-1$. Because when $-1\leq b<0$ and $a > b-1$, a fixed point may have two negative eigenvalues or complex eigenvalues which do not give an invariant half line going to infinity, so we can't apply Brouwer's translation theorem. Below, we concentrate on the case $b>0$. \\

\begin{figure}[htbp]
\centerline{
\includegraphics[width=130mm]{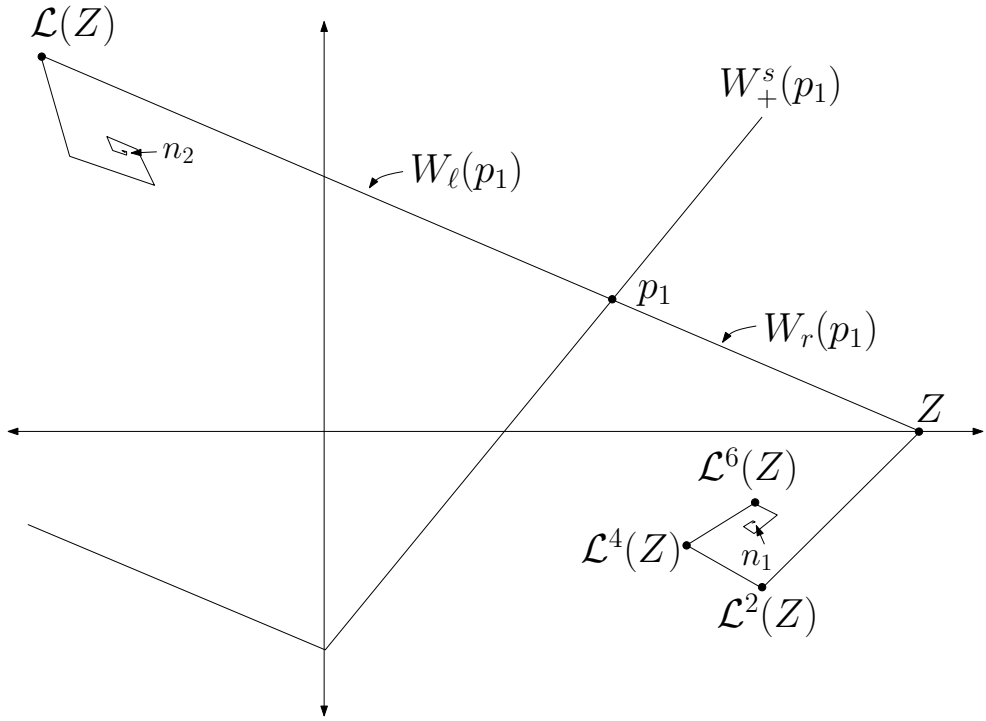}
}
\caption{The picture shows the unstable and stable manifolds of the right fixed point of $\mathcal L_{1,0.5}$.}
\label{unstable picture}
\end{figure}

Let us start stating our results by the following theorem:
\begin{Theorem} For $a=1$ and $b=0.5$, $h_{top}(\mathcal L_{a,b})=0$.
\end{Theorem}
\begin{proof} First note that when $0<b<1$ and $1-b<a<b+1$, $\mathcal L_{a,b}$ has two saddle fixed points: $p_1=(1/(1+a-b),1/(1+a-b))$ in the first quadrant and $p_2=(1/(1-a-b),1/(1-a-b))$ in the third quadrant. Also, there are two attracting period-two points: $n_1=(N, (1-aN)/(1-b))$ in the fourth quadrant and $n_2=((1-aN)/(1-b), N)$ in the second quadrant where $N=(1+a-b)/[(b-1)^2+a^2]$. By a direct calculation of $\mathcal L_{a,b}^4$, one can check that there are no other period-four points. \\

Now $v^{s}_1=(\lambda^{s}_1, 1)$ where $\lambda^{s}_1=(-a+\sqrt{a^2+4b})/2$ is the stable direction at $p_1$ and $W^{s}_{+}(p_1)=\{p_1+v^{s}_{1}t \in \mathbb{R}^2 | \thinspace t>0\}$ is invariant under $\mathcal L_{a,b}$. Similarly, $v^{u}_2=(-\lambda^{u}_2, -1)$ where $\lambda^{u}_2=(a+\sqrt{a^2+4b})/2$ is the unstable direction at $p_2$ and $W^{u}_{+}(p_2)=\{p_2+v^{u}_{2}t \in \mathbb{R}^2 | \thinspace t>0\}$ is invariant under $\mathcal L_{a,b}$. \\

The more challenging part is to show that the right and left parts of the unstable manifold of $p_1$ are attracted by $n_1$ and $n_2$, respectively. We will show that this happens by considering $\mathcal L^4$. Now, let $Z$ be the intersection of the line $\ell_1=\{p_1+v^{u}_1 t \in\mathbb{R}^2 | \thinspace t>0 \}$ and the $x$-axis where $v^{u}_1=(-\lambda^{u}_1, -1)$ and $\lambda^{u}_1=(-a-\sqrt{a^2+4b})/2$. See Fig.~\ref{unstable picture}. \\

\noindent \emph{Claim: For $a=1$ and $b=0.5$, $\mathcal L_{a,b}^{4m}(Z)\rightarrow n_1$ as $m\rightarrow \infty$.} \\

\noindent \emph{Proof of the claim:} Let us use $\mathcal L_{1,0.5}=\mathcal L$. Let $P$ be the polygon whose corners are given by $Z$, $\mathcal L^2(Z)$, $\mathcal L^4(Z)$ and $\mathcal L^6(Z)$. Since $\mathcal L^8(Z)\approx(1.223,-0.375)$ is in $P$, $\mathcal L^2(P)\subset P$, i.e., $P$ is invariant under $\mathcal L^2$. Now consider the Lyapunov function $V(x,y)=(x-\pi_1(n_1))^2+(y-\pi_2(n_1))^2$ where $\pi_1:\mathbb{R}^2 \to \mathbb{R}$ and $\pi_2:\mathbb{R}^2 \to \mathbb{R}$ are the projections to the $x$-coordinate and $y$-coordinate, respectively. Note that $\pi_1(n_1)=6/5$ and $\pi_2(n_1)=-2/5$. By direct calculation one can see that $V(\mathcal L^4(x,y))-V(x,y)= -\dfrac{15}{16}((x-6/5)^2+(y+2/5)^2)<0 \thickspace, \forall (x,y) \in P\setminus\{n_1\}$. This implies that (see for ex. \cite{Saber}), $Z$ (actually every $(x,y) \in P\setminus\{n_1\}$) is asymptotically stable to $n_1$ under $\mathcal L^4$. \\

Similarly it can be shown that $\mathcal L(Z)$ is asymptotically stable to $n_2$ under iterations of $\mathcal L^4$. Now let $W_r(p_1)$ be the union of forward iterations (under $\mathcal L^4$) of the line segment connecting $p_1$ and $Z$, i.e., $W_r(p_1)=\cup_{n\geq0}\mathcal L^{4n}(segment \thinspace[p_1 Z])$. Similarly let $W_\ell (p_1)$ be the union of forward iterations (under $\mathcal L^4$) of the line segment connecting $p_1$ and $\mathcal L(Z)$. To complete the proof of the theorem, we apply Brouwer's translation theorem to $\mathcal L^4$. Note that $\mathbb{R}^2\setminus(W^s_{+}(p_1)\cup\{p_1\}\cup W^u_{+}(p_2)\cup\{p_2\}\cup W_r(p_1)\cup\{n_1\}\cup W_\ell(p_1)\cup\{n_2\})$ is homeomorphic to $\mathbb{R}^2$ and $\mathcal L^4$ has no fixed points there. Since $\mathcal L^4$ is orientation preserving $h_{top}(\mathcal L)=h_{top}(\mathcal L^4)/4=0$.
\end{proof}

\emph{Proof of the Theorem \ref{zero entropy}:} The proof of the above theorem, using similar Lyapunov functions, works for the parameters in a small neighborhood of $(a,b)=(1,0.5)$ as well. \\

\begin{figure}[htbp]
\centerline{
\includegraphics[width=100mm]{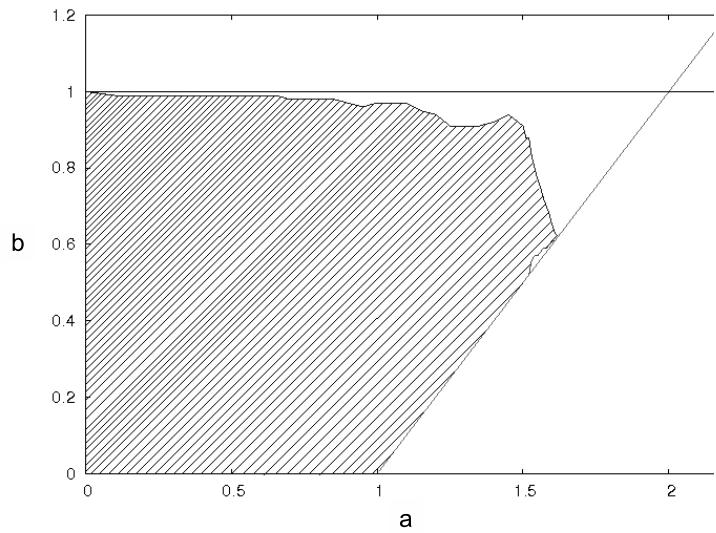}
}
\caption{The shaded region gives the numerically observed parameters for which $h_{top}(\mathcal L_{a,b})=0$.}
\label{zero parameters}
\end{figure}

\emph{Remark:} When we move away from a neighborhood of $(a,b)=(1,0.5)$, it is sometimes the case that the unstable manifold of the right fixed point intersects with the stable manifold of the same fixed point causing a homoclinic point and positive entropy. The parameters for which $\mathcal L_{a,b}$ is numerically observed to have zero entropy are given in Fig.~\ref{loziparameterspace1} and Fig.~\ref{zero parameters}. For more details see \cite{BurakWeb}. Note that since positive entropy occurs as a result of a homoclinic intersection of the stable and unstable manifolds of a periodic point (which are piecewise linear), the boundary of the zero entropy locus is expected to be piecewise algebraic. But writing the equations explicitly requires more work.

\paragraph{The case \textbf{a=1+b}:} When $a=1+b$ and $b>0$, it can be shown that the portion of the line $\ell:\thinspace y=-x+(1-b^2)/(a(1+b^2))$ that stays in the region given by $1+ax+by\geq0$, $1-a(1+ax+by)+bx\leq0$, $x\leq0$ and image of that portion of the line $\ell$ under $\mathcal L_{a,b}$ give all the period-four points except the fixed points of $\mathcal L_{a,b}$. In other words, there are infinitely many period-four points that lie on two line segments. But it can be again observed numerically that as long as there are no homoclinic points, the unstable manifold of the right fixed point is attracted by these two line segments causing the entropy to be zero. Note that when $a>1+b$, the period-two points become saddles, so we can expect that some portion of the line $a=1+b$, $b>0$ is a part of the boundary of the zero entropy locus (see Fig.~\ref{zero parameters}). For more detailed study of this case, see \cite{Burak2}. \\

Finally, we would like to mention some entropy results for the H\'enon family. In \cite{Arai}, Arai introduced a rigorous computational method to show the uniform hyperbolicity of the H\'enon family at several regions in the parameter space. This means the entropy is constant at those regions. Recently Frongillo \cite{Frongillo}, using the techniques introduced in \cite{Day}, gave rigorous lower bounds for the entropy of the H\'enon family at those parameters studied by Arai. These studies gave some insight about the approximate maximal entropy and zero entropy regions for the H\'enon family. Especially, one should compare the conjectural maximal entropy parameters in the H\'enon family \cite{Frongillo} with the maximal entropy parameters in the Lozi family (Fig.~\ref{loziparameterspace1}).
\clearpage

\newpage
\appendix
\section{Appendix}

\begin{Proposition}\label{q} For $\ueps$ with $\uepss=(+1,-1,-1,-1\dots)$ we have $q(\ueps)(a,b)=\frac{b}{(a+x)(b+x)}$ where $x=(a-\sqrt{a^2+4b})/2$. \\
\end{Proposition}
\begin{proof} First note that since $\uepss=(+1,-1,-1,-1\dots)$, \\

  $r_0=\cfrac{1}{a + \cfrac{b}{-a + \cfrac{b}{-a + \cfrac{b}{\ddots}}}} = \cfrac{1}{a + x}$ where $x=\cfrac{b}{-a + \cfrac{b}{-a + \cfrac{b}{\ddots}}}$. \\
  Note that the continued fraction for $x$ can be written as $x=\cfrac{b}{-a + x}$ and this equation gives two solutions. We choose $x=(a-\sqrt{a^2+4b})/2$ as in \cite{Ish1}.

  Also note that $r_n=\cfrac{1}{-a + \cfrac{b}{-a + \cfrac{b}{-a + \cfrac{b}{\ddots}}}}$ for all $n\geq1$. So $r_n=\frac{x}{b}$ for $n\geq1$.

  Now we have $q(\ueps)(a,b)=r_0-r_0 r_1+r_0 r_1 r_2-\ldots= \frac{1}{a+x}(1-\frac{x}{b}+\frac{x^2}{b^2}-\ldots)=\frac{1}{a+x}(\sum\limits_{n=0}^{\infty}(-\frac{x}{b})^n)=\frac{b}{(a+x)(b+x)}$.
\end{proof}

\begin{Lemma}\label{delq} For $1<a\leq2$ and $\ueps \in \sPh_{2,0}$, $\eval[1]{\frac{\partial q}{\partial b}}_{b=0}=\lim\limits_{b \to 0}{\frac{\partial q}{\partial b}}=\lim\limits_{b \to 0}{\frac{\partial }{\partial b}}\left(\frac{b}{(a+x)(b+x)}\right)=\frac{\left(1-\frac{2}{a}\right)}{a(a-1)^2}$. For $a=2$, we have $\lim\limits_{b \to 0}{\frac{\partial q}{\partial b}}=0$.
\end{Lemma}
\begin{proof} Note that by the above proposition, we have $\eval[1]{\frac{\partial q}{\partial b}}_{b=0}=\lim\limits_{b \to 0}{\frac{\partial q}{\partial b}}=\lim\limits_{b \to 0}{\frac{\partial }{\partial b}}\left(\frac{b}{(a+x)(b+x)}\right)$.

Note that $\frac{\partial q}{\partial b}=\frac{(a+x)(b+x)-b[x'(b+x) + (a+x)(1+x')]}{(a+x)^2(b+x)^2}$. To find $\lim\limits_{b \to 0}{\frac{\partial q}{\partial b}}$, we need to apply l'H\^opital's Rule twice. Applying l'H\^opital's Rule the first time, $\lim\limits_{b \to 0}{\frac{\partial q}{\partial b}}$, after cancelation, becomes:

\[ \lim\limits_{b \to 0}\frac{-b[x''(b+x)+x'(1+x')+x'(1+x')+(a+x)x'']}{2(a+x)x'(b+x)^2+(a+x)^22(b+x)(1+x')},
\]

which equals:

\[\lim\limits_{b \to 0}\frac{-b[x''(a+b+2x)+2x'(1+x')]}{2(a+x)(b+x)[(x'(b+x)+(a+x)(1+x')]}.
\]

Applying l'H\^opital's Rule again, the $b$-derivative of the numerator becomes:

\[-[x''(a+b+2x)+2x'(1+x')]-b[x'''(a+b+2x)+x''(1+2x')+2x''(1+x')+2x'x''],
\]

and the $b$-derivative of the denominator becomes:

\begin{eqnarray}
2[(x'(b+x)+(a+x)(1+x'))(x'(b+x)+(a+x)(1+x')) \nonumber\\
  +(a+x)(b+x)(x''(b+x)+2x'(1+x')+(a+x)x'')]. \nonumber
\end{eqnarray}

Now, taking the limit of the numerator and denominator as $b$ goes to $0$ gives the result. Note that $\lim\limits_{b \to 0}{x} = \lim\limits_{b \to 0}\frac{a-\sqrt{a^2+4b}}{2} = 0$ and $\eval[1]{\frac{\partial x}{\partial b}}_{b=0} = \eval[1]{x'}_{b=0}=-\frac{1}{a}$ and $\eval[1]{\frac{\partial^2 x}{\partial b^2}}_{b=0}=\eval[1]{x''}_{b=0}=\frac{2}{a^3}$:

\[ \lim\limits_{b \to 0}{\frac{\partial q}{\partial b}}=\frac{\frac{2}{a}\left(1-\frac{2}{a}\right)}{2a^2\left(1-\frac{1}{a}\right)^2}=\frac{\left(1-\frac{2}{a}\right)}{a(a-1)^2}.
\]
\end{proof}

\paragraph{Acknowledgments}
I would like to thank S.E.~Newhouse for his helpful discussions and suggestions and especially for sharing his computer code that made the results in Section~\ref{results2} possible.

\bibliography{references}
\bibliographystyle{plain}
\ttfamily
Department of Mathematics, Michigan State University, East Lansing, MI 48824 USA \\
Current address: Max Planck Institute for Cognitive and Brain Sciences, Leipzig, Germany \\ E-mail: yildiz@cbs.mpg.de
\end{document}

%% file: macros.tex
\newcommand{\ind}{\mbox{\mbox{ } \mbox{ } \mbox{ }}}

\newcommand{\be}{\begin{enumerate}}
\newcommand{\ee}{\end{enumerate}}
\newcommand{\bit}{\begin{itemize}}
\newcommand{\eit}{\end{itemize}}
\newcommand{\bc}{\begin{center}}
\newcommand{\ec}{\end{center}}
\newcommand{\beq}{\begin{equation}}
\newcommand{\eeq}{\end{equation}}
\newcommand{\bqn}{\begin{eqnarray}}
\newcommand{\eqn}{\end{eqnarray}}
\newcommand{\bqns}{\begin{eqnarray*}}
\newcommand{\eqns}{\end{eqnarray*}}
\newcommand{\N}{{\bf N}}

\newcommand{\R}{{\bf R}}

\newcommand{\rarrow}{\rightarrow}

\newcommand{\sL}{{\mathcal L}}

\newcommand{\se}{{\mathcal E}}

\newcommand{\bdoc}{\begin{document}}
\newcommand{\edoc}{\end{document}}
\newcommand{\docart}{\documentstyle[12pt]{article}}
\newcommand{\vs}{\vspace{.16in}}
\newcommand{\half}{\frac{1}{2}}
\newcommand{\third}{\frac{1}{3}}
\newcommand{\fourth}{\frac{1}{4}}
\newcommand{\diffeo}{\mbox{ diffeomorphism }}
\newcommand{\dnorm}[1]{\mbox{$\mid \mid #1 \mid \mid$}}
\newcommand{\zero}{ \mbox{$ {\bf 0 } $}}
\newcommand{\PL}{\mbox{$Per_L(\R)$}}

\newcommand{\0}{{\bf 0}}

\newcommand{\twoover}[2]{{\footnotesize \begin{array}[t]{c}\mbox{ {\normalsize #1}}  \\
                               #2
        \end{array} } \hspace{1cm} }
\newcommand{\threeover}[3]{ {\footnotesize \begin{array}[t]{c}\mbox{ {\normalsize #1}}  \\
                               #2  \\
                                 #3
        \end{array} } \hspace{1cm} }
\newcommand{\fourover}[4]{ {\footnotesize \begin{array}[t]{c}\mbox{ {\normalsize #1}}  \\
                               #2  \\
                                 #3 \\
				#4
        \end{array} } \hspace{1cm}  }

\newenvironment{vover3}[3]{
       {\footnotesize \begin{array}[t]{c}
                               {\normalsize #1} \vspace{.2cm} \\
                                 #2 \\
                                 #3
                       \end{array} }
                            } { \hspace{1cm} }
\newenvironment{vover4}[4] {
       {\footnotesize \begin{array}[t]{c}
                               {\normalsize #1} \vspace{.2cm} \\
                                 #2 \\
                                 #3 \\
                                 #4
                       \end{array} }
                            } {  \hspace{1cm} }

\newenvironment{seteqover}[3]{#1 \ \ = \hspace{.6cm}
                                #2 #3} {\ind}

\newcommand{\Norm}[1]{\mbox{$\parallel #1 \parallel$}}
\newcommand{\imp}{\mbox{$\Longrightarrow$}}
\newcommand{\pr}{\mbox{$\prime$}}

\newcommand{\twocvec}[2]{\mbox{$\left( \begin{array}{c}
				 #1 \\
				 #2
				\end{array} \right)$} }

\newcommand{\twomat}[4]{\mbox{$ \left( \begin{array}{cc}
                              #1 & #2 \\
                           #3 & #4
                       \end{array} \right) $ } }

\newcommand{\Itwomat}[4]{\mbox{$ \left( \begin{array}{cc}
                              #1 & #2 \\
                           #3 & #4
                       \end{array} \right)^{-1} $ } }

\newcommand{\threemat}[9]{ \[ \left( \begin{array}{ccc}
				 #1 & #2 & #3 \\
				 #4 & #5 & #6 \\
				 #7 & #8 & #9
				\end{array} \right) \]  }

\newcommand{\newthreemat}[9]{ \left( \begin{array}{ccc}
				 #1 & #2 & #3 \\
				 #4 & #5 & #6 \\
				 #7 & #8 & #9
				\end{array} \right) }

\newcommand{\Ithreemat}[9]{\mbox{ $\left( \begin{array}{ccc}
				 #1 & #2 & #3 \\
				 #4 & #5 & #6 \\
				 #7 & #8 & #9
				\end{array} \right)^{-1} $ } }

\newcommand{\ba}{\mbox{${\bf a}$}}
\newcommand{\bfa}{\mbox{${\bf a}$}}
\newcommand{\bb}{\mbox{${\bf b}$}}
\newcommand{\bfb}{\mbox{${\bf b}$}}
\newcommand{\bfc}{\mbox{${\bf c}$}}
\newcommand{\bi}{\mbox{${\bf i}$}}
\newcommand{\bj}{\mbox{${\bf j}$}}
\newcommand{\bfi}{\mbox{${\bf i}$}}
\newcommand{\bfj}{\mbox{${\bf j}$}}
\newcommand{\bfk}{\mbox{${\bf k}$}}
\newcommand{\bfx}{\mbox{${\bf x}$}}
\newcommand{\bfy}{\mbox{${\bf y}$}}
\newcommand{\bfz}{\mbox{${\bf z}$}}
\newcommand{\bfu}{\mbox{${\bf u}$}}
\newcommand{\bfv}{\mbox{${\bf v}$}}
\newcommand{\bfw}{\mbox{${\bf w}$}}
\newcommand{\bfg}{\mbox{${\bf g}$}}
\newcommand{\onef}{\mbox{${\bf 1}$}}
\renewcommand{\bfc}{\mbox{${\bf c}$}}
\newcommand{\bfe}{\mbox{${\bf e}$}}
\newcommand{\bff}{\mbox{${\bf f}$}}
\newcommand{\bS}{\mbox{${\bf S}$}}
\newcommand{\bV}{\mbox{${\bf V}$}}
\newcommand{\bw}{\bar{w}}
\newcommand{\fh}{\hat{f}}
\newcommand{\wh}{\hat{w}}
\newcommand{\Wh}{\hat{W}}
\newcommand{\Ih}{\hat{I}}
\newcommand{\Lh}{\hat{L}}
\newcommand{\grgh}{\hat{\gamma}}
\newcommand{\br}{\mbox{${\bf r}$}}
\newcommand{\bs}{\mbox{${\bf s}$}}

\newcommand{\oz}{\overline{z}}
\newcommand{\ogz}{\overline{gz}}
\newcommand{\oxi}{\overline{\xi}}

\newcommand{\ninf}{n \rarrow \infty}
\newcommand{\linf}{\lim_{n \rarrow \infty}}
\newcommand{\fon}{ \frac{1}{n} }
\newcommand{\skn}{ \sum_{k=0}^{n-1} }
\newcommand{\siin}{ \sum_{i=0}^{n-1} }

\newcommand{\alp}{\mbox{$\alpha$}}
\newcommand{\gra}{\mbox{$\alpha$}}

\newcommand{\bet}{\mbox{$\beta$}}
\newcommand{\grb}{\mbox{$\beta$}}

\newcommand{\gam}{\mbox{$\gamma$}}
\newcommand{\Gam}{\mbox{$\Gamma$}}
\newcommand{\Del}{\mbox{$\Delta$}}
\newcommand{\sig}{\mbox{$\sigma$}}
\newcommand{\Sig}{\mbox{$\Sigma$}}

\newcommand{\grg}{\mbox{$\gamma$}}
\newcommand{\Grg}{\mbox{$\Gamma$}}
\newcommand{\grd}{\mbox{$\delta$}}
\newcommand{\Grd}{\mbox{$\Delta$}}
\newcommand{\grs}{\mbox{$\sigma$}}
\newcommand{\Grs}{\mbox{$\Sigma$}}
\newcommand{\grf}{\mbox{$\phi$}}
\newcommand{\Grf}{\mbox{$\Phi$}}
\newcommand{\grth}{\mbox{$\theta$}}
\newcommand{\Grth}{\mbox{$\Theta$}}
\newcommand{\gro}{\mbox{$\omega$}}
\newcommand{\Gro}{\mbox{$\Omega$}}
\newcommand{\gre}{\mbox{$\epsilon$}}
\newcommand{\bgre}{\mbox{$\bar{\epsilon}$}}
\newcommand{\grve}{\mbox{$\varepsilon$}}
\newcommand{\grr}{\mbox{$\rho$}}
\newcommand{\Grr}{\mbox{$\Rho$}}
\newcommand{\grt}{\mbox{$\tau$}}
\newcommand{\gri}{\mbox{$\iota$}}
\newcommand{\grp}{\mbox{$\pi$}}
\newcommand{\grP}{\mbox{$\Pi$}}
\newcommand{\grl}{\mbox{$\lambda$}}
\newcommand{\grlh}{\mbox{$\hat{\lambda}$}}
\newcommand{\grL}{\mbox{$\Lambda$}}
\newcommand{\Grl}{\mbox{$\Lambda$}}
\newcommand{\grz}{\mbox{$\zeta$}}

\newcommand{\thet}{\mbox{$\theta$}}
\newcommand{\Thet}{\mbox{$\Theta$}}
\newcommand{\ome}{\mbox{$\omega$}}
\newcommand{\Ome}{\mbox{$\Omega$}}
\newcommand{\eps}{\mbox{$\epsilon$}}

\newcommand{\iot}{\mbox{$\iota$}}

\newcommand{\lam}{\mbox{$\lambda$}}
\newcommand{\Lam}{\mbox{$\Lambda$}}
\newcommand{\zet}{\mbox{$\zeta$}}

\newcommand{\kap}{\mbox{$\kappa$}}

\newcommand{\grx}{\mbox{$\xi$}}
\newcommand{\Grx}{\mbox{$\Xi$}}
\newcommand{\grc}{\mbox{$\chi$}}
\newcommand{\Grc}{\mbox{$\Chi$}}
\newcommand{\grn}{\mbox{$\nu$}}
\newcommand{\grm}{\mbox{$\mu$}}
\newcommand{\grk}{\mbox{$\kappa$}}

\newtheorem{theorem}{Theorem}[section]
\newtheorem{Cor}[theorem]{Corollary}
\newtheorem{Corollary}[theorem]{Corollary}
\newtheorem{Prop}[theorem]{Proposition}
\newtheorem{Proposition}[theorem]{Proposition}
\newtheorem{Definition}[theorem]{Definition}
\newtheorem{Lem}[theorem]{Lemma}
\newtheorem{Lemma}[theorem]{Lemma}
\newtheorem{Rem}[theorem]{Remark}
\newtheorem{Remark}[theorem]{Remark}
\newtheorem{Conjecture}[theorem]{Conjecture}
\newtheorem{Theorem}[theorem]{Theorem}
\newtheorem{Th}[theorem]{Theorem}
\newtheorem{question}{Question}
\newtheorem{conjecture}{Conjecture}

\newcommand{\ta}{\tilde{a}}
\newcommand{\tgra}{\tilde{\alpha}}
\newcommand{\tGrl}{\tilde{\Grl}}
\newcommand{\tb}{\tilde{b}}
\newcommand{\tgrb}{\tilde{\beta}}
\newcommand{\tgrg}{\tilde{\gamma}}
\newcommand{\tc}{\tilde{c}}
\newcommand{\td}{\tilde{d}}
\newcommand{\te}{\tilde{e}}
\newcommand{\tf}{\tilde{f}}
\newcommand{\teta}{\tilde{\eta}}
\newcommand{\tg}{\tilde{g}}
\newcommand{\tih}{\tilde{h}}
\renewcommand{\th}{\tilde{h}}
\newcommand{\ti}{\tilde{i}}
\newcommand{\tj}{\tilde{j}}
\newcommand{\tk}{\tilde{k}}
\newcommand{\tl}{\tilde{l}}
\newcommand{\tm}{\tilde{m}}
\newcommand{\tn}{\tilde{n}}
\newcommand{\tio}{\tilde{o}}
\newcommand{\tp}{\tilde{p}}
\newcommand{\tq}{\tilde{q}}
\newcommand{\tr}{\tilde{r}}
\newcommand{\ts}{\tilde{s}}
\newcommand{\tit}{\tilde{t}}
\newcommand{\tu}{\tilde{u}}
\newcommand{\tv}{\tilde{v}}
\newcommand{\tw}{\tilde{w}}
\newcommand{\tx}{\tilde{x}}
\newcommand{\ty}{\tilde{y}}
\newcommand{\tz}{\tilde{z}}
\newcommand{\tA}{\tilde{A}}
\newcommand{\tB}{\tilde{B}}
\newcommand{\tC}{\tilde{C}}
\newcommand{\tD}{\tilde{D}}
\newcommand{\tE}{\tilde{E}}
\newcommand{\tF}{\tilde{F}}
\newcommand{\tG}{\tilde{G}}
\newcommand{\tH}{\tilde{H}}
\newcommand{\tI}{\tilde{I}}
\newcommand{\tJ}{\tilde{J}}
\newcommand{\tK}{\tilde{K}}
\newcommand{\tL}{\tilde{L}}
\newcommand{\tM}{\tilde{M}}
\newcommand{\tN}{\tilde{N}}
\newcommand{\tO}{\tilde{O}}
\newcommand{\tP}{\tilde{P}}
\newcommand{\tQ}{\tilde{Q}}
\newcommand{\tR}{\tilde{R}}
\newcommand{\tS}{\tilde{S}}
\newcommand{\tT}{\tilde{T}}
\newcommand{\tU}{\tilde{U}}
\newcommand{\tV}{\tilde{V}}
\newcommand{\tW}{\tilde{W}}
\newcommand{\tX}{\tilde{X}}
\newcommand{\tY}{\tilde{Y}}
\newcommand{\tZ}{\tilde{Z}}

\newcommand{\tTh}{\tilde{\Theta}}

\newcommand{\fl}{f_{\grl}}
\newcommand{\el}{E_{\grl}}
\newcommand{\sel}{\se_{\grl}}
\newcommand{\szl}{\sz_{\grl}}
\newcommand{\zl}{Z_{\grl}}
\newcommand{\famfl}{\{ \fl \} }
\newcommand{\sT}{{\mathcal T}}
\newcommand{\hsT}{{\hat{\mathcal T}}}
\newcommand{\hsH}{{\hat{\mathcal H}}}
\newcommand{\hsCZ}{{\hat{\mathcal CZ}}}
\newcommand{\sCG}{{\mathcal CG}}
\newcommand{\sFG}{{\mathcal FG}}
\newcommand{\sTemp}{{\mathcal TM}}
\newcommand{\ab}{{\bf a}}
\newcommand{\xb}{{\bf x}}
\newcommand{\yb}{{\bf y}}
\newcommand{\Def}{{\bf Definition}}
\newcommand{\Pf}{{\bf Proof}}
\newcommand{\uG}{\bar{{\mathcal G}}}
\newcommand{\sG}{{\mathcal G}}
\newcommand{\sF}{{\mathcal F}}
\newcommand{\stG}{\tilde{{\mathcal G}}}
\newcommand{\tmu}{\tilde{\mu}}
\newcommand{\tpi}{\tilde{\pi}}
\newcommand{\tgrl}{\tilde{\grl}}
\newcommand{\sD}{{\mathcal D}}
\newcommand{\sDh}{\hat{\sD}}
\newcommand{\sP}{{\mathcal P}}
\newcommand{\tsP}{{\tilde{\mathcal P}}}
\newcommand{\sPh}{\hat{\sP}}
\newcommand{\sM}{{\mathcal M}}
\newcommand{\bsM}{\bar{{\mathcal M}}}
\newcommand{\sN}{{\mathcal N}}
\newcommand{\sS}{{\mathcal S}}
\newcommand{\tsS}{\tilde{{\mathcal S}}}
\newcommand{\tgrt}{\tilde{{\tau}}}
\newcommand{\sI}{{\mathcal I}}
\newcommand{\sW}{{\mathcal W}}
\newcommand{\sCZ}{{\mathcal C}{\mathcal Z}}
\newcommand{\tsT}{\tilde{\sT}}
\newcommand{\tsC}{\tilde{\sC}}
\newcommand{\tsD}{\tilde{\sD}}
\newcommand{\tgre}{\tilde{\gre}}
\newcommand{\sH}{{\mathcal H}}
\newcommand{\tsH}{\tilde{\sH}}
\newcommand{\tsG}{\tilde{\sG}}
\newcommand{\tsW}{\tilde{\sW}}
\newcommand{\tDel}{\tilde{\Grd}}
\newcommand{\sA}{{\mathcal A}}
\newcommand{\sB}{{\mathcal B}}
\newcommand{\sCF}{{\mathcal CF}}
\newcommand{\bsH}{\bar{\sH}}
\newcommand{\ntwob}{\bar{\nu_2}}
\newcommand{\nb}{\bar{\nu}}
\newcommand{\bE}{\bar{E}}
\newcommand{\bF}{\bar{F}}
\newcommand{\bG}{\bar{G}}
\newcommand{\bI}{\bar{I}}
\newcommand{\bZ}{\bar{Z}}
\newcommand{\bh}{\bar{h}}
\newcommand{\bm}{\bar{m}}
\newcommand{\ml}{m_{\ell}}
\newcommand{\bT}{\bar{T}}
\newcommand{\bK}{\bar{K}}
\newcommand{\bx}{\bar{x}}
\newcommand{\bB}{\bar{B}}
\newcommand{\bxi}{\bar{\xi}}
\newcommand{\bu}{\bar{u}}
\newcommand{\bv}{\bar{v}}
\newcommand{\bg}{\bar{g}}
\newcommand{\bN}{\bar{\N}}
\newcommand{\bX}{\bar{X}}
\newcommand{\bY}{\bar{Y}}
\newcommand{\bGrs}{\bar{\Grs}}
\newcommand{\bgrs}{\bar{\grs}}
\newcommand{\bmu}{\bar{\mu}}
\newcommand{\bnu}{\bar{\nu}}
\newcommand{\bTh}{\bar{\Theta}}
\newcommand{\db}{\bar{d}}
\newcommand{\dx}{\dot{x}}
\newcommand{\dy}{\dot{y}}
\newcommand{\dz}{\dot{z}}
\newcommand{\ddx}{\partial_x}
\newcommand{\ddy}{\partial_y}
\newcommand{\ddz}{\partial_z}
\newcommand{\dxn}{\dot{x}_n}
\newcommand{\n}{\norm}
\newcommand{\No}{\dnorm}
\newcommand{\no}[1]{\mbox{$\norm{#1}_1$}}
\newcommand{\grfh}{\hat{\grf}}

\newcommand{\ol}{\overline}
\newcommand{\ob}{\bar}

\newcommand{\eqdef}{\mbox{ $\stackrel{{\rm def}}{=}$ }}
\newcommand{\defeq}{\mbox{ $\stackrel{{\rm def}}{=}$ }}

\newcommand{\ueps}{\underline{\grve}}
\newcommand{\uepss}{\underline{\grve}^s}
\newcommand{\uepsu}{\underline{\grve}^u}
\newcommand{\headeps}{(\varepsilon_{0},\varepsilon_{1} \ldots )}
\newcommand{\taileps}{(\ldots, \varepsilon_{-2},\varepsilon_{-1})}
\newcommand{\fulleps}{(\ldots \varepsilon_{-2},\varepsilon_{-1}\cdot\varepsilon_{0},\varepsilon_{1}\ldots) }

\newcommand{\lex}[2]{{\vphantom{#2}}_{#1}{#2}} 

\newcommand{\Lpfk}{\mbox{$\cap \hspace{-.44cm}\mid$}}